\documentclass[11pt]{amsart}

\usepackage{amsmath,amssymb,amsthm,amsfonts,amscd,flafter,epsf,epsfig,graphicx,verbatim}
\usepackage[all]{xy}

\title[Genus one, one boundary component open books]{Heegaard Floer homology and genus one, one boundary component open books}
\author[John A. Baldwin]{John A. Baldwin}
\address {Department of Mathematics, Columbia University\\ New York, NY 10027}
\email {baldwin@math.columbia.edu}
\date{}

\newcommand\ww{\mathbb{W}}
\newcommand\tp{\mathcal{T}^+}
\newcommand\ms{M_{T,\phi}} 
\newcommand\mt{\widehat{M}_{T,\phi}}

\newcommand\Sc{\text{Spin}^c}
\newcommand\spc{\mathfrak{s}}
\newcommand\spt{\mathfrak{t}}
\newcommand\so{\mathfrak{s}_0}
\newcommand\kh{\widetilde{Kh}}

\newcommand\hfi{HF^{\infty}}

\newcommand\hfp{HF^{+}}
\newcommand\hf{\widehat{HF}}

\newcommand\cf{\widehat{CF}}

\newcommand\cfki{CFK^{\infty}}
\newcommand\hfk{\widehat{HFK}}
\newcommand\zz{\mathbb{Z}}
\newcommand\zt{\mathbb{F}}
\newcommand\zzt{\mathbb{Z}_2}
\newcommand\Sig{S}
\newcommand\Tight{\text{Tight}(T, \partial T)}
\newcommand\Dehn{\text{Dehn}^+(T, \partial T)}

\newcommand\Stein{\text{Stein}(T, \partial T)}
\newcommand\al{\mathbf{\alpha}}
\newcommand\be{\mathbf{\beta}}

\newcommand\x{\sigma_1}
\newcommand\y{\sigma_2}
\newcommand\Kp{K_+}
\newcommand\Kv{K_v}
\newcommand\Kh{K_h}

\newcommand\uknot{\mathcal{U}}
\newcommand\tref{\mathcal{T}}
\newcommand\fige{\mathcal{E}}

\newtheorem{theorem}{Theorem}[section]
\newtheorem{lemma}[theorem]{Lemma}

\newtheorem{corollary}[theorem]{Corollary}
\newtheorem{proposition}[theorem]{Proposition}

\theoremstyle{definition}
\newtheorem{definition}[theorem]{Definition}

\newtheorem{remark}[theorem]{Remark}

\begin{document}
\begin{abstract}  
We compute the Heegaard Floer homology of any rational homology 3-sphere with an open book decomposition of the form $(T,\phi)$, where $T$ is a genus one surface with one boundary component. In addition, we compute the Heegaard Floer homology of every $T^2$-bundle over $S^1$ with first Betti number equal to one, and we compare our results with those of Lebow on the embedded contact homology of such torus bundles. We use these computations to place restrictions on Stein-fillings of the contact structures compatible such open books, to narrow down somewhat the class of 3-braid knots with finite concordance order, and to identify all quasi-alternating links with braid index at most 3.  \end{abstract} 

\maketitle

\section{Introduction}
An open book is a pair $(\Sig, \phi)$, where $\Sig$ is a compact surface with boundary, and $\phi$ is a diffeomorphism of $\Sig$ which restricts to the identity on $\partial \Sig$. From an open book, one may construct a closed, oriented 3-manifold $M_{\Sig,\phi}$ as follows: $M_{\Sig,\phi} = \Sig \times [0,1] /\sim$, where $\sim$ is the identification defined by 
\begin{eqnarray*}
(x,1) \sim(\phi(x),0), && x \in \Sig\\
(x,t) \sim (x, s), && x\in \partial \Sig, \ t,s \in [0,1].
\end{eqnarray*}

We say that $(\Sig,\phi)$ is an \emph{open book decomposition} for the 3-manifold $Y$ if $Y\cong M_{\Sig,\phi}$. The image of $\partial \Sig \times \{t\}$ in $M_{\Sig,\phi}$ is called the \emph{binding} of this open book decomposition, and the image of $S \times \{t\}$ is called a \emph{page}.  The binding is a fibered link in $M_{\Sig,\phi}$ with genus $g = g(\Sig)$. Conversely, every genus $g$ fibered link $K \subset Y$ is isotopic to the binding of some open book decomposition $(\Sig,\phi)$ of $Y$, with $g(\Sig)=g$. 

By performing $0$-surgery on each component of the binding, we obtain an $\widehat{\Sig}$-bundle over $S^1$, where $\widehat{\Sig}$ is the closed surface obtained by capping off each boundary component of $\Sig$ with a disk. This bundle has monodromy $\widehat{\phi}$, where $\widehat{\phi}$ is the extension of $\phi$ to $\widehat{\Sig}$ by the identity. Conversely, every $\widehat{\Sig}$-bundle over $S^1$ may be obtained in this way. We denote this fibered 3-manifold by $\widehat{M}_{\Sig,\phi}$.

In this paper, we study the Heegaard Floer homology of the 3-manifolds $\ms$, where $T$ is a genus one surface with one boundary component. Specifically, we compute the $\zz[U]$-module $\hfp(\ms)$ for all $\phi$ for which $b_1(\ms)=0$. These 3-manifolds are precisely the rational homology 3-spheres which contain genus one fibered knots, or, equivalently, those 3-manifolds which arise as branched double covers of $S^3$ along closed 3-braids with non-zero determinant. The Heegaard Floer homology of $\ms$ is most interesting in the torsion $\Sc$ structure, $\so$, associated to the contact structure $\xi_{T,\phi}$ compatible with $(T,\phi)$ (in the sense of Giroux \cite{giroux}). For this $\Sc$ structure, we determine the absolute $\mathbb{Q}$-grading on $\hfp(\ms,\so)$. For all other $\spc \in \Sc(\ms)$, we show that $\hfp(\ms,\spc)$ is isomorphic to $\hfp(S^3) \cong \zz[U,U^{-1}]/(U\cdot\zz[U])$ as a relatively graded $\zz[U]$-module.

Along the way, we compute the Heegaard Floer homology of the torus bundles $\mt$ with $b_1(\mt)=1$. Torus bundles are of interest in other Floer theories as well. For instance, when $\phi$ is the identity, Hutchings and Sullivan show in \cite{hs} that the embedded contact homology of $\mt \cong T^3$ agrees with its Heegaard and Seiberg-Witten Floer homologies up to shifts in grading, evidence for the conjecture that these three Floer theories are equivalent despite their very different constructions. Eli Lebow has since computed the embedded contact homology for almost all $T^2$-bundles over $S^1$ \cite{lebow}. In Section \ref{section:tb.msmt}, we compare our results with his and verify that embedded contact homology and Heegaard Floer homology are isomorphic as relatively graded $\zz$-modules for any $T^2$-bundle over $S^1$ with \emph{pseudo-Anosov} monodromy. 

Surface bundles over $S^1$ have also been studied in the context of Heegaard Floer homology by Jabuka and Mark in \cite{jabmark, jabmark2, jabmark3}, and by Roberts in \cite{lrob2}. As we shall see, our computations may be used to extend some of their results on the Heegaard Floer homology in non-torsion $\Sc$ structures of $\Sigma_g$-bundles over $S^1$, where $\Sigma_g$ is a closed surface of genus $g>1$.

With our computations of Heegaard Floer homology, we identify all \emph{$L$-spaces} among the manifolds $\ms$. 

\begin{definition} An $L$-space is a rational homology 3-sphere whose Heegaard Floer homology is as simple as possible; namely, $\hf(Y,\spc) \cong \zz$, or, equivalently, $\hfp(Y,\spc) \cong \zz[U,U^{-1}]/(U\cdot\zz[U]),$ for all $\spc \in \text{Spin}^c(Y)$. 
\end{definition} There is no purely topological characterization of $L$-spaces, though $L$-spaces do come with some rigid \emph{geometric} restrictions. Most notably, Ozsv{\'a}th and Szab{\'o} prove that if $Y$ is an $L$-space then $Y$ contains no co-orientable taut foliation \cite{osz2}. In \cite{bald1}, we use this result to produce an infinite family of hyperbolic 3-manifolds with no co-orientable taut foliations, extending the first known set of such 3-manifolds, which was discovered by Roberts, Shareshian, and Stein among fillings of punctured torus bundles \cite{rss}. 

Let $\Tight$ (respectively, $\Stein$) denote the set of boundary-fixing diffeomorphisms $\phi$ of $T$ for which $\xi_{T,\phi}$ is tight (respectively, Stein-fillable), and let us denote by $\Dehn$ the set of diffeomorphisms which are isotopic to a product of right-handed Dehn twists around non-separating curves in $T$. In this paper, we use the $\mathbb{Q}$-grading on $\hf(\ms,\so)$ to provide a formula for the Euler characteristic of any Stein-filling of $\xi_{T,\phi}$ whenever $\ms$ is an $L$-space. We use this formula to find an infinite family of diffeomorphisms $\phi \in \Tight-\Stein$. This improves upon work of Honda, Kazez, and Mati{\'c} who, in \cite{hkm2}, identify an infinite family of diffeomorphisms $\phi \in \Tight-\Dehn$.

\begin{figure}[!htbp]
\begin{center}
\includegraphics[width=9cm]{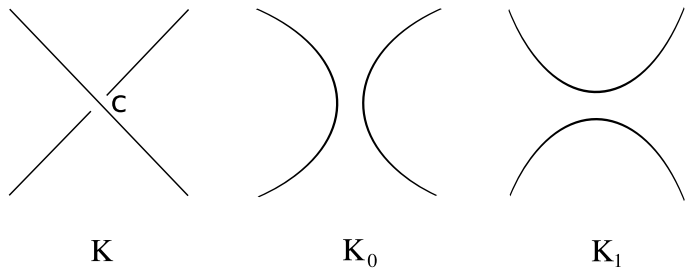}
\caption{The resolutions $K_0$ and $K_1$ at the crossing $c$.}
\label{fig:Res}
\end{center}
\end{figure}

$L$-spaces are commonly found among branched covers of $S^3$. Ozsv{\'a}th and Szab{\'o} prove, for instance, that if $K \subset S^3$ is an alternating link then $\Sigma(K)$, the double cover of $S^3$ branched along $K$, is an $L$-space \cite{osz12}. In the same paper, they define a more general class of links called \emph{quasi-alternating} links whose branched double covers are also $L$-spaces:

\begin{definition} The set $Q$ of quasi-alternating links is the smallest set of links which satisfies the following properties: 
\label{def:Q}
\begin{itemize}
\item The unknot is in $Q$,
\item If $K$ admits a projection with a crossing $c$ for which 
\begin{enumerate}
\item both resolutions, $K_0, K_1$, as in Figure \ref{fig:Res}, are in $Q$, and
\item $\text{det}(K) = \text{det}(K_0) + \text{det} (K_1),$
\end{enumerate} then $K\in Q$.
\end{itemize}
\end{definition} All alternating links are quasi-alternating. Yet, quasi-alternating knots abound among small crossing non-alternating knots \cite{man1}. 

Quasi-alternating links widen the bridge between knot Floer homology and Khovanov homology. Specifically, Manolescu and Ozsv{\'a}th show in \cite{manozs} that if $K$ is quasi-alternating then both its reduced Khovanov homology, $\kh(K)$, and its knot Floer homology, $\hfk(S^3,K)$, are \emph{$\sigma$-thin}; that is, both theories are supported in bi-gradings $(i,j)$ which satisfy $j-i = -\sigma(K)/2$, where $\sigma(K)$ is the signature of $K$. This was first shown to be the case for alternating knots in \cite{osz9, lee}. It follows that if $K$ is a quasi-alternating $l$-component link, then \begin{eqnarray*}\text{rk}(\kh(K)) = \text{det}(K)& \text{     and     } &\text{rk}(\hfk(S^3,K)) = 2^{l-1}\text{det}(K).\end{eqnarray*} 

Ultimately, as in the case of $L$-spaces, there is no purely topological description of quasi-alternating links. In Section \ref{section:tb.3braid}, we give a complete description of 3-braids whose closures are quasi-alternating using, in part, our classification of $L$-spaces among the manifolds $\ms$. In particular, we prove:

\begin{proposition}
\label{prop:thin}
If $K$ is a link with braid index at most 3 then $K$ is quasi-alternating if and only if its reduced Khovanov homology $\kh(K)$ is $\sigma$-thin.
\end{proposition}

It seems likely that there are links $K$ which are not quasi-alternating but for which $\kh(K)$ is $\sigma$-thin. No such examples have been confirmed in print, but the knot $9_{46}$ is a good candidate.\footnote{The author has recently learned that Shumakovitch has shown that the reduced odd Khovanov homology (see \cite{orsz}) of $9_{46}$ is not thin, and, hence, that $9_{46}$ is not quasi-alternating.}

Our computations of $\mathbb{Q}$-graded Heegaard Floer homology can also be used to elicit information about the smooth concordance order of 3-braid knots. 
\begin{definition}
If $c_1(\spc)$ is torsion then the correction term of $Y$ in the $\Sc$ structure $\spc$ is defined to be the grading of the lowest degree element in the image of the map $\pi:\hfi(Y,\spc) \rightarrow \hfp(Y,\spc)$. It is denoted by $d(Y,\spc)$.
\end{definition} When $K$ is a knot, $\Sigma(K)$ supports a unique self-conjugate $\Sc$ structure $\so$. In \cite{mo}, Manolescu and Owens show that the correction term $d(\Sigma(K),\so)$ provides information about the smooth concordance order of $K$. We compute these correction terms for all 3-braid knots (since, in this case, $\Sigma(K)$ is diffeomorphic to some $\ms$ with $b_1(\ms)=0$), and use the Manolescu-Owens machinery along with the knot signature to narrow down the class of 3-braid knots with finite concordance order in Section \ref{section:tb.3braid}.

\subsection*{Organization} 
In Section \ref{section:mcg}, we explain the correspondence between 3-braids and genus one, one boundary component open books. In Section \ref{section:tb.surg}, we describe the Heegaard Floer homology of $p/q$-surgery on a genus one fibered knot in an $L$-space. In Section \ref{section:tb.lspace}, we classify all $L$-spaces among the manifolds $\ms$. In Section \ref{section:corr}, we compute correction terms for many of the $L$-spaces found in Section \ref{section:tb.lspace}. In Section \ref{section:tb.msmt}, we give an explicit description of the graded Heegaard Floer homologies of $\ms$ and $\mt$, and we compare this description with Lebow's computations of embedded contract homology. In addition, we extend work of Jabuka and Mark on the Heegaard Floer homology of $\Sigma_g$-bundles over $S^1$. In Section \ref{section:stein}, we study Stein-fillings of the contact structures $\xi_{T,\phi}$. Finally, in Section \ref{section:tb.3braid}, we discuss applications to the smooth concordance order of 3-braid knots, and we classify all quasi-alternating links with braid index at most 3.

A computer program which implements many of the results of this paper is available at \begin{tt}http://www.math.columbia.edu/$\sim$baldwin/3braid.html\end{tt}. 

\subsection*{Acknowledgements}
I wish to thank Joan Birman, Josh Greene, Ciprian Manolescu, and Peter Ozsv{\'a}th for helpful discussions. In addition, I am grateful to Eli Lebow for sharing his unpublished computations of embedded contact homology.
\newpage

\section{Open books of the form $(T,\phi)$ and 3-braids}
\label{section:mcg}

\begin{figure}[!htbp]
\begin{center}
\includegraphics[height=8cm]{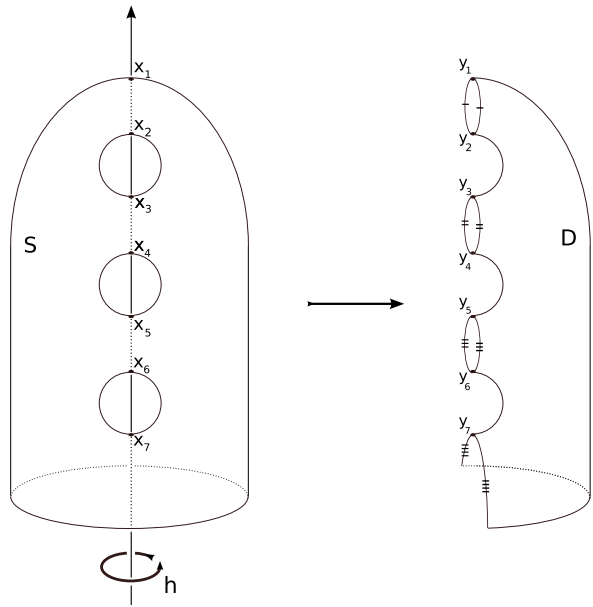}
\caption{On the left is the surface $S$. We've drawn the case in which the genus of $S$ is 3. On the right is the quotient of $S$ by the involution $h$. When the boundary arcs are identified as indicated, $D$ is a disk.}
\label{fig:hyperelliptic}
\end{center}
\end{figure}

Suppose that $S$ is a compact surface with one boundary component, and let $h$ denote the hyperelliptic involution of $S$ obtained from a $180^o$ rotation about the axis depicted in Figure \ref{fig:hyperelliptic}. Let $\text{Aut}(S,\partial S)$ denote the set of isotopy classes of orientation preserving diffeomorphisms of $S$ which fix $\partial S$ pointwise (where isotopies also fix $\partial S$ pointwise). If $\phi \in \text{Aut}(S,\partial S)$ commutes with $h$, then $h$ extends to an involution $\hat{h}$ of $M_{S,\phi}$ via $\hat{h}((x,t)) = (h(x),t)$. The map $\hat{h}$ is well-defined since $$\hat{h}((x,1)) = (h(x),1) = (\phi \circ h(x),0) = (h \circ \phi (x),0) = \hat{h}((\phi(x),0)).$$ The quotient of $M_{S,\phi}$ by $\hat{h}$ has an open book decomposition $(D, \psi)$, where $D$ is the quotient of $S$ by $h$, and $\psi$ is the diffeomorphism of $D$ induced by $\phi$. From Figure \ref{fig:hyperelliptic} it is clear that $D$ is a disk; hence, $(D, \psi)$ is an open book decomposition for $S^3$, and, therefore, $M_{S,\phi}$ is a branched double cover of $S^3$. 

The branch locus in $S^3$ consists of those points $p$ with only one pre-image in $M_{S,\phi}$ under the covering map $M_{S,\phi} \rightarrow M_{S,\phi}/\hat{h} \cong S^3$. Using the induced open book decomposition $(D,\psi)$ of $S^3$, we see that this branch locus consists precisely of the images of intervals $\{y\}\times[0,1]$ for which $y \in D$ has only one pre-image in $S$ under the covering map $S \rightarrow S/h \cong D$. But $y\in D$ has only one pre-image $x\in S$ if and only if $h(x) = x$. When the genus of $S$ is $g$ then $h(x) = x$ precisely for the $2g+1$ points $x_1, \dots, x_{2g+1}$ which are at the intersection of $S$ with the rotational axis depicted in Figure \ref{fig:hyperelliptic}. Hence, if we let $y_i$ denote the image of $x_i$ under the covering map, as illustrated in Figure \ref{fig:hyperelliptic}, then the branch locus in $S^3$ consists precisely of the images of the intervals $\{y_i\}\times[0,1]$ in the open book decomposition $(D,\psi)$ for $S^3$. Therefore, this branch locus is the closure of a braid with $2g+1$ strands. The braid is specified by $\psi$ and, hence, by $\phi$.

\begin{figure}[!htbp]
\begin{center}
\includegraphics[height=4cm]{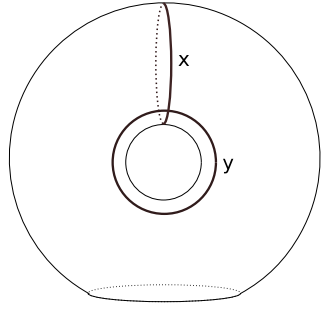}
\label{fig:MCG}
\caption{The surface $T$ and the dual, non-separating curves, $x$ and $y$.}
\end{center}
\end{figure}

Since $\text{Aut}(T,\partial T)$ is generated by Dehn twists around the curves $x$ and $y$ depicted in Figure \ref{fig:MCG}, it is clear that every diffeomorphism of the surface $T$ is isotopic to a diffeomorphism which commutes with the involution $h$. Therefore, the 3-manifolds $\ms$ are exactly those which arise as branched double covers of $S^3$ along closed 3-braids. The lemma below expresses this correspondence between 3-braids and open books $(T,\phi)$ more concretely.

\begin{lemma}
\label{lem:DBC}
Let $\x$ and $\y$ denote the elementary generators of $B_3$, the braid group on three strands, and suppose that $w=w(\x,\y,\x^{-1},\y^{-1})$ is a word in these generators and their inverses. If $K$ is the closed 3-braid specified by $w$ then $\Sigma(K)$ has an open book decomposition $(T, w(x,y,x^{-1},y^{-1}))$.
\end{lemma} 

The following classification of 3-braids is due to Murasugi \cite{mur}.

\begin{theorem}
\label{thm:BraidClass}
Let $w$ be a braid word in $B_3$, and let $h = (\x\y)^3$. Then $w$ is equivalent in $B_3$, after conjugation, to one of the following.
\begin{enumerate}
\item $h^d \cdot \x\y^{-a_1} \cdots \x\y^{-a_n}$, where the $a_i \geq 0$, some $a_j \neq 0$.
\item $h^d \cdot \y^m$, for $m \in \mathbb{Z}.$
\item $h^d \cdot \x^m \y^{-1},$ where $m \in \{-1,-2,-3\}.$
\end{enumerate}
Note that $h$ is a full positive twist of the braid. 
\end{theorem} 

In addition to Lemma \ref{lem:DBC}, it is helpful to know that $\text{Aut}(T,\partial T)$ is isomorphic to $B_3$ via the homomorphism which sends $x$ to $\x$ and $y$ to $\y$. The corollary below follows immediately.

\begin{corollary}
\label{cor:GOFClass}
Let $h = (xy)^3$. Every 3-manifold with a genus one, one boundary component open book decomposition is diffeomorphic to $\ms$ where $\phi$ is one of the following.
\begin{enumerate}
\item $h^d \cdot xy^{-a_1} \cdots xy^{-a_n}$, where the $a_i \geq 0$, some $a_j \neq 0$.
\item $h^d \cdot y^m$, for $m \in \mathbb{Z}.$
\item $h^d \cdot x^m y^{-1},$ where $m \in \{-1,-2,-3\}.$
\end{enumerate}
\end{corollary}

An independent proof of this fact can be found in \cite{bald1}. There, we give an algorithm which takes any word in the Dehn twists $x$ and $y$ and their inverses, and converts it into one of the above forms. 
\newpage

\section{Surgery on genus one fibered knots in $L$-spaces}
\label{section:tb.surg}

We shall see in Section \ref{section:tb.lspace} that every 3-manifold $\ms$ with $b_1(\ms) =0$ is obtained via $1/n$-surgery on a genus one fibered knot in an $L$-space. Accordingly, this section is devoted to understanding the effect of such surgeries on Heegaard Floer homology. More generally, we give a preliminary description of the Heegaard Floer homology of any 3-manifold obtained via $p/q$-surgery on such a knot. When $p/q \neq 0$ these include all manifolds with $b_1=0$ obtained by Dehn-filling a once-punctured torus bundle. When $p/q=0$, these include all $T^2$-bundles over $S^1$ with $b_1=1$. This section relies on the interaction of contact geometry and Heegaard Floer homology. Many details are omitted here, but can be found in \cite{osz1} and \cite{et2}. 

\subsection{Knot Floer homology and the higher differentials}
\label{ssection:tb.surg.hfk}
Unless otherwise specified, we will work with coefficients in a field $\zt$ throughout this section. 

A co-orientable contact structure on a 3-manifold $Y$ is a nowhere integrable 2-plane field $\xi$ which arises as the kernel of a 1-form $\alpha$ with $d\alpha \land \alpha >0$. In \cite{giroux}, Giroux proves that there is a 1-1 correspondence between isotopy classes of contact structures on $Y$ and open book decompositions of $Y$ up to an equivalence called positive stabilization. Using this correspondence, Ozsv{\'a}th and Szab{\'o} associate to a contact structure $\xi$ on $Y$ an element $c(\xi) \in \hf(-Y)/\pm 1$ which is an invariant of $\xi$ up to isotopy. This invariant is constructed as follows. 

Recall that a knot $K$ in a 3-manifold $Y$ induces a filtration on $\cf(Y)$. The filtered chain homotopy type of this complex is an invariant of the knot $K$, and $\hfk(Y,K)$ is the homology of the associated graded object arising from this filtration. Therefore, there is a spectral sequence associated to $K \subset Y$ whose $E^1$ term is $\hfk(Y,K)$, and which converges to $\hf(Y)$. 

Let $(\Sig, \phi)$ be an open book decomposition of $Y$ with connected binding which is compatible with $\xi$. As described in the introduction, this binding $K \subset Y \cong M_{\Sig,\phi}$ is a fibered knot with genus $g = g(\Sig)$. Reversing orientation, we can think of $K$ as a knot in $-M_{\Sig,\phi}$. According to Ozsv{\'a}th and Szab{\'o}, the knot Floer homology of a genus $g$ fibered knot has rank 1 in extremal Alexander gradings $\pm g$: 
\begin{equation}
\label{eqn:Fib}
\hfk(-M_{\Sig,\phi},K, g) \cong \hfk(-M_{\Sig,\phi},K, -g) \cong \zt.
\end{equation}  The generator of $\hfk(-M_{\Sig,\phi},K, -g)$ represents a cycle throughout the spectral sequence associated to $K \subset -M_{\Sig,\phi}$, and its image in $\hf(-M_{\Sig,\phi})$ is the contact invariant $c(\xi)$, which we also denote by $c(\Sig,\phi)$ (it is well-defined up to sign). With this in place, we may begin.

\begin{proposition}
\label{prop:HFK}
Suppose that $\ms$ is an $L$-space, and let $K$ denote the binding of the open book decomposition $(T,\phi)$. Furthermore, let $\so$ denote the $\Sc$ structure associated to the contact 2-plane field corresponding to $(T,\phi)$. If $\spc\neq \so$ then
$$ \hfk(M_{T,\phi},K,i,\spc) \cong 
\left\{\begin{array}{ll}
\zt, &\text{if $i=0,$}\\
0, &\text{otherwise.}
\end{array}\right.$$ If $\spc=\so$ then $\hfk(M_{T,\phi}, K, \pm 1, \so) \cong \zt$. Below, we have listed the three possibilities for $\hfk(M_{T,\phi}, K, 0, \so)$ and the $d^1$ differentials.
 
\begin{enumerate}
\item If $c(T,\phi) \neq 0$ then $\hfk(M_{T,\phi}, K, 0, \so) \cong \zt$, and there is one non-trivial differential $d^1: \hfk(M_{T,\phi}, K, 0, \so) \rightarrow \hfk(M_{T,\phi}, K, -1, \so).$ \\

\item If $c(T,\phi^{-1}) \neq 0$ then $\hfk(M_{T,\phi}, K, 0, \so) \cong \zt$, and there is one non-trivial differential $d^1: \hfk(M_{T,\phi}, K, 1, \so) \rightarrow \hfk(M_{T,\phi}, K, 0, \so).$\\

\item If $c(T,\phi)=c(T,\phi^{-1}) = 0$ then $\hfk(M_{T,\phi}, K, 0, \so) \cong \zt^3$, and there are two non-trivial differentials, $d_j^1: \hfk(M_{T,\phi}, K, j, \so) \rightarrow \hfk(M_{T,\phi}, K, j-1, \so),$ for $j=0,\, 1$.
 \end{enumerate}
 The spectral sequence from $\hfk(\ms,K,\so)$ to $\hf(\ms,\so)$ collapses at the $E^2$ term in each case.
\end{proposition} 

\begin{remark}
\label{remark:ss} Let $\uknot$, $\tref$, $-\tref$, and $\fige$ denote the unknot, the right-handed trefoil, the left-handed trefoil, and the figure-eight, respectively. Proposition \ref{prop:HFK} says that in any $\Sc$ structure $\spc \neq \so$ the $(E^1,d^1)$ term of the spectral sequence for $K \subset M_{T,\phi}$ is identical to that of the spectral sequence for $\uknot\subset S^3$. Meanwhile, in the $\Sc$ structure $\so$, the $(E^1,d^1)$ term of the spectral sequence for $K$ is identical to that of the spectral sequence for $\tref$ when $c(T,\phi) \neq 0$, $-\tref$ when $c(T,\phi^{-1})$, and $\fige$ when $c(T,\phi) = c(T,\phi^{-1})=0$.

\end{remark}

Moreover, when $c(T,\phi) \neq 0$ or $c(T,\phi^{-1}) \neq 0$, the relative Maslov grading on the knot Floer homology of $K$ is immediately clear from the conjugation symmetry discussed below, combined with the fact that all spectral sequence differentials lower Maslov grading by 1. In particular, this relative grading is the same as the relative Maslov grading on the knot Floer homology of $\tref$ or $-\tref$, respectively. When $c(T,\phi)=c(T,\phi^{-1})=0$, the relative grading is the same as the relative grading on the knot Floer homology of $\fige$, but this is less obvious. Nonetheless, this fact may be gleaned from the proof of Proposition \ref{prop:Correction}.(2).

\begin{proof}[Proof of Proposition \ref{prop:HFK}] Since $K$ is a genus one fibered knot in $-\ms$, Equation \ref{eqn:Fib} implies that $\hfk(-\ms,K,-1) \cong \zt$. Ozsv{\'a}th and Szab{\'o} show that the generator of $\hfk(-\ms,K,-1)$ is supported in the $\Sc$ structure $\so$. Moreover, Etnyre and Ozbagci prove in \cite{et3} that this $\Sc$ structure $\so$ is self-conjugate.\footnote{In fact, whenever $\nu$ is a diffeomorphism of a genus $g$ surface $\Sig$ with one boundary component, which commutes with the hyperelliptic involution (as depicted in Figure \ref{fig:hyperelliptic}) on $\Sig$, the $\Sc$ structure associated to the contact structure $\xi_{\Sig,\nu}$ is self-conjugate \cite{et3}. } By the conjugation symmetry in knot Floer homology \cite{osz3}, $$\hfk(-M_{T,\phi}, K, -1, \so) \cong \hfk(-M_{T,\phi}, K, 1, \bar{\spc}_0) \cong \hfk(-M_{T,\phi}, K, 1, \so),$$ the generator of $\hfk(-M_{T,\phi}, K, 1)$ is also supported in the $\Sc$ structure $\so$. 

Now, there is a filtered chain complex associated to $K \subset M_{T,\phi}$ which is dual to the complex associated to $K \subset -M_{T,\phi}$ (obtained, for instance, by switching the orientation of the Heegaard surface). It follows that $\hfk(\ms,K,\pm 1,\so)\cong \zt$. Thus, for $\spc\neq \so$, $\hfk(\ms,K,\pm 1,\spc)=0,$ which implies that the spectral sequence whose $E^1$ term is $\hfk(\ms,K,\spc)$, and which converges to $\hf(\ms,\spc) \cong \zt$ (since $\ms$ is an $L$-space) collapses at the $E^1$ term when $\spc\neq \so$. That is, $\hfk(\ms,K,0,\spc) \cong \zt.$ 

In the $\Sc$ structure $\so$, this spectral sequence collapses at the $E^2$ term. For, the higher differentials $d^i$ lower Maslov grading by 1 and Alexander grading by $i$. Since $\hfk(\ms,K,\so)$ is supported in Alexander gradings -1, 0, and 1, the differentials $d^i$ vanish for $i>2$. Moreover, due to the conjugation symmetry $\hfk_d(\ms,K,-i,\so) \cong \hfk_{d-2i}(\ms,K,i,\so),$ the $d^2$ differentials must be trivial as well, as the generators of $\hfk(\ms,K,\pm 1,\so)$ have Maslov gradings which differ by 2.

Since this spectral sequence collapses at the $E^2$ term to $\hf(\ms,\so) \cong \zt$, it follows from Euler characteristic considerations that $\hfk(\ms,K,0,\so) \cong \zt^m$ where $m=1$ or 3. Below we analyze $\hfk(\ms,K,0,\so)$ and the $d^1$ differentials for the three cases in Proposition \ref{prop:HFK}. \\

\noindent \textbf{I.} Suppose $c(T,\phi) \neq 0$. Then the generator of $\hfk(-\ms,K,-1,\so)$ survives in the spectral sequence associated to $K \subset -\ms$, and represents the generator of the $E^2$ term. Thus, there is exactly one non-trivial differential, $$d^1: \hfk(-M_{T,\phi}, K, 1, \so) \rightarrow \hfk(-M_{T,\phi}, K, 0, \so).$$ By the duality discussed above, there is one non-trivial differential, $$d^1: \hfk(M_{T,\phi}, K, 0, \so) \rightarrow \hfk(M_{T,\phi}, K, -1, \so).$$ Hence, $\hfk(\ms,K,0,\so)\cong \zt$ since $H_*(E^1,d^1)$ must be isomorphic to $\zt$.\\

\noindent \textbf{II.} Suppose $c(T,\phi^{-1}) \neq 0$. Observe that $(T,\phi^{-1})$ is an open book decomposition for the 3-manifold $M_{T,\phi^{-1}} \cong -M_{T,\phi}$. Since $c(T,\phi^{-1}) \neq 0$, there is exactly one non-trivial differential, $$d^1: \hfk(M_{T,\phi}, K, 1, \so) \rightarrow \hfk(M_{T,\phi}, K, 0, \so).$$ Hence, $\hfk(\ms,K,0,\so)\cong \zt$. \\

\noindent \textbf{III.} Suppose $c(T,\phi)=0$ and $c(T,\phi^{-1})=0$. Since $c(T,\phi)=0$, the generator of the group $\hfk(-\ms,K,-1,\so)$ is killed by a non-trivial differential, $$d^1: \hfk(-M_{T,\phi}, K, 0, \so) \rightarrow \hfk(-M_{T,\phi}, K, -1, \so).$$ By duality, there is a non-trivial differential, $$d^1: \hfk(M_{T,\phi}, K, 1, \so) \rightarrow \hfk(M_{T,\phi}, K, 0, \so).$$  Again, $(T,\phi^{-1})$ is an open book decomposition for $-\ms$, and since $c(T,\phi^{-1})=0$, there is a non-trivial differential, $$d^1: \hfk(M_{T,\phi}, K, 0, \so) \rightarrow \hfk(M_{T,\phi}, K, -1, \so).$$ Hence, $\hfk(\ms,K,0,\so) \cong \zt^3.$

\end{proof}

\subsection{Surgery}
\label{ssection:tb.surg.surg}

This subsection is dedicated to an understanding of the $\mathbb{Q}$-graded $\zz[U]$-module structure of $\hfp((\ms)_{p/q}(K))$ when $\ms$ is an $L$-space and $K$ is the binding of the open book $(T,\phi)$. First we establish a bit of notation.

If $G$ is a graded group, we let $G\{n\}$ denote the group obtained from $G$ by shifting all gradings by $n$. Let $\tp_{d}$ denote the graded $\zz[U]$-module $\zz[U,U^{-1}]/(U\cdot \zz[U])$ in which the element $1$ has grading $d$, and multiplication by $U$ drops grading by 2. Observe that we can identify $\Sc((\ms)_{p/q}(K))$ with $\Sc(\ms) \oplus \zz/p\zz$. Therefore, we denote a $\Sc$ structure on $(\ms)_{p/q}(K)$ by $(\spc,i)$, where $\spc \in \Sc(\ms)$ and $i \in \zz / p\zz$.  Below is the main result of this section.

\begin{proposition}
\label{prop:Surgery}
Suppose that $\ms$ is an $L$-space, and let $K$ denote the binding of the open book decomposition $(T,\phi)$. When $\spc =\so$, $$\hfp((\ms)_{p/q}(K), (\so,i)) \cong
\left\{\begin{array}{ll} \hfp(S^3_{p/q}(\tref),i)\{d(\ms,\so)\} \\ 
\hfp(S^3_{p/q}(-\tref),i)\{d(\ms,\so)\}  \\ 
\hfp(S^3_{p/q}(\fige),i)\{d(\ms,\so)\}  
 \end{array} \right.$$
 The three lines on the right-hand side above correspond to the three cases $c(T,\phi)\neq 0$, $c(T,\phi)=c(T,\phi)= 0$, and $c(T,\phi^{-1})\neq 0$, respectively. For all $\spc \neq \so$, $$\hfp((\ms)_{p/q}(K),(\spc,i)) \cong \hfp(S^3_{p/q}(\uknot),i)\{d(\ms, \spc)\}.$$
\end{proposition}

\begin{proof}[Proof of Proposition \ref{prop:Surgery}] 
Recall that the chain complex $\cfki(Y,K, \spc)$ is a free $\zt[U,U^{-1}]$-module generated by tuples of the form $[\mathbf{x},i,j]$, where $\mathbf{x} \in \mathbb{T}_{\alpha} \cap \mathbb{T}_\beta$ corresponds to the $\Sc$ structure $\spc $, and the Alexander grading of $\mathbf{x}$ (thought of as an element of the filtered complex $\cf(Y,\spc)$) is $A(\mathbf{x}) = j-i$. In fact, by restricting to a vertical (or horizontal) slice of $\cfki(Y,K, \spc)$, we recover the filtered complex $\cf(Y, \spc)$, where the Alexander grading of a generator is given by $j-i$ (or $i-j$).

In \cite{osz7}, Ozsv{\'a}th and Szab{\'o} prove that for any $p/q$, $\hfp(Y_{p/q}(K), (\spc,t))$ is determined as a relatively graded $\zt[U]$-module by the bi-filtered chain homotopy type of $\cfki(Y,K, \spc)$. Proposition \ref{prop:Surgery} follows from the fact that the filtered chain homotopy type of this complex is determined, for the knots that we are considering, by the $(E^1,d^1)$ term of the spectral sequence for $K\subset Y$. 

Let $\widehat{\partial}_{i}$ be the component of the differential $\widehat{\partial}$ on the complex $\cf(Y,K, \spc)$ which drops the filtration grading by $i$. The ``cancellation lemma" (see Rasmussen's thesis \cite{ras}) ensures that there exists a complex $(\widehat{C}(Y, \spc),\widehat{\partial}')$ which is filtered chain homotopic to $(\cf(Y, \spc),\widehat{\partial})$, and for which $$\widehat{C}(Y, \spc) \cong H_*(\cf(Y, \spc), \widehat{\partial}_{0}).$$ That is, $$\widehat{C}(Y, \spc) \cong \hfk(Y,K, \spc)$$ 

We can apply the cancellation lemma to the bi-filtered complex $\cfki(Y,K, \spc)$ in the same way. Indeed, let $\partial_{ij}$ be the component of the differential $\partial$ on $\cfki(Y,K, \spc)$ which drops the bi-grading by $(i,j)$. Then there exists a complex $(C^{\infty}(Y,K, \spc),\partial')$ which is bi-filtered chain homotopic to $(\cfki(Y,K, \spc),\partial)$, and for which $$C^{\infty}(Y,K, \spc) \cong H_*(\cfki(Y,K, \spc), \partial_{0,0})$$ as an $\zt[U,U^{-1}]$-module. 

By construction, a vertical (or horizontal) slice of the complex $(C^{\infty}(Y,K, \spc),\partial')$ is naturally isomorphic to the complex $(\widehat{C}(Y, \spc),\widehat{\partial}')$, and the differential $\partial'_{0j}$ (or $\partial'_{j0}$) restricted to this slice is simply $\widehat{\partial}'_{j}$. In particular, there is an isomorphism of relatively graded $\zt[U,U^{-1}]$-modules, \begin{equation}\label{eqn:HFK} C^{\infty}(Y,K, \spc) \cong \hfk(Y,K,\spc)\otimes \zt[U,U^{-1}].\end{equation} 

If $K\subset \ms$ is a knot as in the hypothesis of Proposition \ref{prop:Surgery} then Equation \ref{eqn:HFK} together with the relative Maslov grading on $\hfk(\ms,K,\spc)$ implies that the differentials $\partial'_{ij}$ lower Maslov grading by $i+j$, and, hence, must vanish whenever $i+j>1$. Therefore, the only non-vanishing components of $\partial'$ are $\partial'_{10}$ and $\partial'_{01}.$ As mentioned above, these restrict to the differential $\widehat{\partial}'_1$ on the complex $(\widehat{C}(\ms, \spc),\widehat{\partial}')$ when looking at vertical or horizontal slices. But the complex $(\widehat{C}(\ms, \spc),\widehat{\partial}'=\widehat{\partial}'_1)$ is isomorphic to the $(E^1,d^1)$ term of the spectral sequence associated to $K\subset \ms$ in the $\Sc$ structure $\spc $, which we computed in Proposition \ref{prop:HFK} (in particular, $\widehat{\partial}'_1=0$ when $\spc\neq \so$). Therefore, this $(E^1,d^1)$ term completely determines the complex $(C^{\infty}(\ms,K, \spc),\partial')$ for the knots in Proposition \ref{prop:Surgery}. 

We illustrate what is left of the proof of Proposition \ref{prop:Surgery} for the case in which $\spc =\so$ and $c(T,\phi)=0$ (the other cases are treated identically). See Figure \ref{fig:T} for an depiction of $(C^{\infty}(\ms,K,\so),\partial')$ in this case, and note that this chain complex is isomorphic to $(\cfki(S^3,\tref),\partial)$ for a specific choice of doubly-pointed Heegaard diagram for $\tref\subset S^3$ \cite{osz3}. Therefore, $$\hfp((\ms)_{p/q}(K),(\so,i)) \cong \hfp(S^3_{p/q}(\tref),i)\{d(\ms,\so)\}.$$ This grading shift comes from further inspection of the Ozsv{\'a}th-Szab{\'o} formula for rational surgery on knots (see \cite{osz7}), combined with the fact that $d(S^3)=0$. 

Since the above isomorphism holds for coefficients in every field $\zt$, the Universal Coefficient Theorem implies that the above groups are isomorphic as graded $\zz$-modules. By the naturality of the Universal Coefficient Theorem, this is actually a $\zz[U]$-module isomorphism. 

\begin{figure}[!htbp]
\begin{center}
\includegraphics[height=6cm]{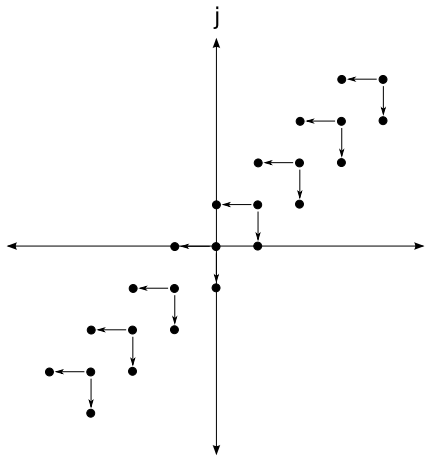}
\caption{The complex $(C^{\infty}(\ms,K,\so),\partial')$ when $c(T,\phi) \neq 0.$ The dots represent generators of $\zt$ at the corresponding coordinates, and the arrows represent the non-trivial components of the $\partial'$ differential.}
\label{fig:T}
\end{center}
\end{figure}

 \end{proof}

%In the case that $p/q = 0$, $(\ms)_0(K)$ is the torus bundle over $S^1$ with monodromy $\phi$. Let us denote this manifold by $\mt$. $\Sc(\mt)$ can be identified with $\Sc(\ms)\oplus \zz$, where $(s,i)$ is a $\Sc$ structure with $\langle c_1((s,i)), [\widehat{S}] \rangle = 2i$, and $\widehat{S}$ is the torus in $\mt$ obtained by capping off the surface $S$. The adjunction inequality for Heegaard Floer homology implies that $\hfp(\mt,(s,i)) = 0$ for all $i \neq 0$ \cite{osz14}. Below we describe $\hfp(\mt,(s,0))$.

%\begin{prop}
%\label{prop:TorusBundle}
%Suppose that $\ms$ is an $L$-space, and let $K$ denote the binding of the open book decomposition $(T,\phi)$. For all $s \neq \so$, we have the following isomorphism of graded $\zz[U]$-modules, $$\hfp(\mt,(s,0)) \cong \hfp(S^1\times S^2)\{d(\ms,s)\}.$$ When $s=\so$ we have, $$\hfp(\mt, (\so,0)) \cong
%\left\{\begin{array}{ll} \hfp(S^3_{0}(\tref))\{d(\ms,\so)\},  &{\text{if $c(\tref,\phi) \neq 0$}}\\ 
%\hfp(S^3_{0}(-\tref))\{d(\ms,\so)\},  &{\text{if $c(T,\phi) = 0,$ $c(T,\phi) \neq 0$}}\\ 
%\hfp(S^3_{0}(\fige))\{d(\ms,\so)\},  &{\text{if $c(T,\phi) = c(T,\phi^{-1})=0$.}}
 %\end{array} \right.$$
%\end{prop}

%\begin{proof}[Proof of Proposition \ref{prop:TorusBundle}] The proof here is identical to that of the preceding proposition. In \cite{osz7}, Ozsv{\'a}th and Szab{\'o} show that $\hfp(\mt, (\so,0))$ is isomorphic to the mapping cone of $v^+_{\so,0} + h^+_{\so,0}:A^+_{\so,0} \rightarrow B^+_{\so}$.
%\end{proof}

In Section \ref{section:tb.msmt}, we use Proposition \ref{prop:Surgery} to give an explicit description of the $\mathbb{Q}$-graded Floer homology $\hfp(\ms,\so)$ for all $\phi$ in Corollary \ref{cor:GOFClass} with $b_1(\ms)=0$. To this end, we first classify those $\phi$ for which $M_{T,\phi}$ is an $L$-space so that we can apply the proposition. The following corollary is useful in this regard.

\begin{corollary}
\label{cor:Lspace}
Let $h=(xy^3)$, and suppose that $\ms$ is an $L$-space. 

\begin{itemize}
\item If $c(T,\phi) \neq 0$ then $M_{T,h^{2d} \cdot \phi}$ is an $L$-space if and only if $d\in\{0,-1\}$.

\item If $c(T,\phi) = c(T,\phi^{-1}) = 0$ then $M_{T,h^{2d} \cdot \phi}$ is an $L$-space if and only if $d=0$.
\end{itemize}
\end{corollary}

\begin{proof}[Proof of Corollary \ref{cor:Lspace}]
If $\gamma$ is a curve in $T$ then the 3-manifold $M_{T,\gamma \cdot w}$ is obtained from $M_{T,w}$ by performing $-1$-surgery on a knot isotopic to $\gamma$ in a page of the open book decomposition $(T,w)$ ($M_{T,\gamma^{-1} \cdot w}$ is obtained via $+1$-surgery). The diffeomorphism $h^2 = (xy)^6$ is isotopic to a right-handed Dehn twist around a curve $\delta$ parallel to $\partial T$. Thus, if $K$ is the binding of the open book decomposition $(T,\phi)$ then $M_{T,h^{2d}\cdot \phi} \cong M_{T,\delta^{d}\cdot \phi}$ is obtained from $M_{T,\phi}$ by performing $-1/d$ surgery on $K$. If $c(T,\phi) \neq 0$ then Proposition \ref{prop:Surgery} implies that $$\hfp(M_{T,h^{2d}\cdot \phi}) \cong \hfp((M_{T,\phi})_{-1/d}(K)) \cong \hfp(S^3_{-1/d}(\tref)).$$ But $S^3_{-1/d}(\tref)$ is an $L$-space if and only if $d\in\{0,-1\}$ \cite{osz6}. Meanwhile, if $c(T,\phi) = c(T,\phi) =0$, then $$\hfp(M_{T,h^{2d}\cdot \phi}) \cong \hfp(S^3_{-1/d}(\fige)).$$ But $S^3_{-1/d}(\fige)$ is an $L$-space if and only if $d=0$ \cite{osz6}.
\end{proof}

\newpage
\section{$L$-spaces among the manifolds $\ms$}
\label{section:tb.lspace}

In \cite{bald1} we classify $L$-spaces among the 3-manifolds $\ms$ for the diffeomorphisms $\phi$ in Corollary \ref{cor:GOFClass}.(1). As mentioned in the introduction, these are the pseudo-Anosov diffeomorphisms of the once-punctured torus. In that paper, we use the taut foliations constructed by Roberts in \cite{roberts1} and \cite{roberts2} to prove that most of these manifolds are not $L$-spaces. In this section, we extend our $L$-space classification to all three types of diffeomorphisms in Corollary \ref{cor:GOFClass}. In addition, we revisit the $L$-space classification for the pseudo-Anosov diffeomorphisms with an alternate approach that does not rely on the existence of taut foliations; rather, we employ Corollary \ref{cor:Lspace}. 

\begin{theorem}
\label{thm:GOFClass}
The following is a complete classification of $L$-spaces among the 3-manifolds $\ms$. 
\begin{enumerate}
\item If $\phi = h^d \cdot xy^{-a_1} \cdots xy^{-a_n}$, where $a_i \geq 0$ and some $a_j \neq 0$, then $\ms$ is an $L$-space if and only if $d \in\{-1,0,1\}$.
\item If $\phi = h^d \cdot y^m$ then $\ms$ is an $L$-space if and only if $d = \pm 1$.
\item If $\phi = h^d \cdot x^m y^{-1},$ where $m \in \{-1,-2,-3\}$, then $\ms$ is an $L$-space if and only if $d \in \{-1,0,1,2\}$.
\end{enumerate}
\end{theorem}

The theorem below is a combination of results from \cite{bald1} and \cite{hkm3}, and will be helpful in the proof of Theorem \ref{thm:GOFClass} as well as in the description of the Heegaard Floer homologies of $\ms$ and $\mt$ in Section \ref{section:tb.msmt}.

\begin{theorem}
\label{thm:Tight}
The contact structure $\xi_{T,\phi}$ is tight if and only if $c(T,\phi)\neq 0$. Therefore, the following is a complete classification of tight contact structures compatible with open books of the form $(T,\phi)$.
\begin{enumerate}
\item If $\phi = h^d \cdot xy^{-a_1} \cdots xy^{-a_n}$, where $a_i \geq 0$ and some $a_j \neq 0$, then $c(T,\phi)\neq0$ if and only if $d >0$.
\item If $\phi = h^d \cdot y^m$ then $c(T,\phi)\neq0$ if and only if $d>0$ or $d=0$ and $m\geq0$.
\item If $\phi = h^d \cdot x^m y^{-1},$ where $m \in \{-1,-2,-3\}$, then $c(T,\phi)\neq0$ if and only if $d >0$.
\end{enumerate}
\end{theorem}

\begin{proof}[Proof of Theorem \ref{thm:GOFClass}] This proof is broken into three parts.\\

\noindent \textbf{I.} Let us first consider the diffeomorphisms $\phi = h^d \cdot xy^{-a_1} \cdots xy^{-a_n},$ where $a_i \geq 0$ and some $a_j \neq 0$. Let $d=0$. Lemma \ref{lem:DBC} implies that $M_{T,xy^{-a_1} \cdots xy^{-a_n}}$ is the double cover of $S^3$ branched along the closed 3-braid specified by $$\x\y^{-a_1} \cdots \x\y^{-a_n}.$$ This braid closure is an alternating link, and, hence, $M_{T,xy^{-a_1} \cdots xy^{-a_n}}$ is an $L$-space. The following property of Heegaard Floer homology, proved in \cite{osz4}, is very useful for the case $d=1$.

\begin{theorem}
\label{thm:LspaceH1}
Suppose $K \subset Y$ is a knot with meridian $\mu$. For any framing $\lambda$, we let $Y_{\lambda}(K)$ denote the 3-manifold obtained via $\lambda$-framed surgery on $K$. If $Y$ and $Y_{\lambda}(K)$ are both $L$-spaces, and $$|H_1(Y_{\mu+\lambda}(K);\mathbb{Z})| = |H_1(Y;\mathbb{Z})| + |H_1(Y_{\lambda}(K);\mathbb{Z})|,$$ then $Y_{\mu+\lambda}(K)$ is an $L$-space.
\end{theorem}

Let $d=1$. As described in the proof of Corollary \ref{cor:Lspace}, the 3-manifold $M_{T,h \cdot xy^{-a_1} \cdots xy^{-a_n}}$ is obtained from $M_{T,h \cdot xy^{-a_1} \cdots xy^{-a_n+1}}$ by performing $+1$-surgery around a curve isotopic to $y$ supported in a page of the open book decomposition for $M_{T,h \cdot xy^{-a_1} \cdots xy^{-a_n+1}}$.\footnote{When $\phi = h \cdot xy^{-a_1} \cdots xy^{-a_n}$ is conjugate to $h\cdot x^ny^{-1}$ we show directly that $M_{T, \phi}$ is an $L$-space. In every other case, we can assume that either $a_n>1$ or that $a_n=1$ and some other $a_i=1$ to ensure that $h \cdot xy^{-a_1} \cdots xy^{-a_n+1}$ is still a diffeomorphism of the sort in Theorem \ref{thm:GOFClass}.(1).}  Keeping with the notation in \cite{bald1}, we let $-Q(a_1+1,a_2,\dots,a_{n-1},a_{n-1}+1)$ denote the result of $0$-surgery around the same curve. Then the manifolds $$M_{T,h \cdot xy^{-a_1} \cdots xy^{-a_n}}, \, M_{T,h \cdot xy^{-a_1} \cdots xy^{-a_n+1}}, \, \text{and } -Q(a_1+1,a_2,\dots,a_{n-1},a_{n-1}+1)$$ are related by surgeries as in the hypothesis of Theorem \ref{thm:LspaceH1}. Moreover, in \cite{bald1}, we showed the following:

\begin{lemma}\label{lem:H1} For $a_i\geq 0$ and some $a_j \neq 0$,
\begin{eqnarray*}
|H_1(M_{T, h \cdot xy^{-a_1} \cdots xy^{-a_n-1}};\zz)| &=& |H_1(M_{T, h \cdot xy^{-a_1} \cdots xy^{-a_n}};\zz)|\\
&+& |H_1(-Q(a_1+1,a_2,\dots,a_{n-1},a_{n-1}+1);\zz)|.
\end{eqnarray*} 
\end{lemma}

In addition, $-Q(a_1+1,a_2,\dots,a_{n-1},a_{n-1}+1)$ is an $L$-space (it is just the branched double cover of $S^3$ along an alternating link). Therefore, we can induct using Theorem \ref{thm:LspaceH1} to prove that every 3-manifold of the form $M_{T,h \cdot xy^{-a_1} \cdots xy^{-a_n}}$ is an $L$-space. We need only to check the base case -- that $M_{T,h\cdot x^my^{-1}}$ is an $L$-space for all $m>0$ -- but this is easy. Modulo conjugation and the chain relation ($xyx=yxy$) in $\text{Aut}(T,\partial T)$, 
\begin{eqnarray*}
h\cdot x^my^{-1}  &=& y^{-1}\cdot (xy)^3\cdot x^m = y^{-1}\cdot xyxyxy\cdot x^m = y^{-1}\cdot yxyxyx\cdot x^m \\
&=& xyxyx^{m+1}=xxyxx^{m+1} = yx^{m+4}. 
\end{eqnarray*} Therefore, $M_{T,h\cdot x^my^{-1}}\cong M_{T,yx^{m+4}}$ is the double cover of $S^3$ branched along the 3-braid $\y\x^{m+4}$. But this braid is the alternating link $T(2,m+4)$. Hence, $M_{T,h\cdot x^my^{-1}}$ is an $L$-space. In fact, it is just the lens space $-L(4+m,1)$. 

Thus far, we have shown that $M_{T,h^d \cdot xy^{-a_1} \cdots xy^{-a_n}}$ is an $L$-space when $d=0,1$. By Theorem \ref{thm:Tight}, $c(T,h \cdot xy^{-a_1} \cdots xy^{-a_n}) \neq 0$. As a result, Corollary \ref{cor:Lspace} implies that the 3-manifold $M_{T,h^{2d+1}\cdot xy^{-a_1} \cdots xy^{-a_n}}$ is an $L$-space if and only if $d\in \{0,-1\}$. Moreover, Theorem \ref{thm:Tight} implies that $$c(T,xy^{-a_1} \cdots xy^{-a_n}) = c(T,(xy^{-a_1} \cdots xy^{-a_n})^{-1}) = 0.$$ Consequently, Corollary \ref{cor:Lspace} says that $M_{T,h^{2d}\cdot xy^{-a_1} \cdots xy^{-a_n}}$ is an $L$-space if and only if $d=0$. Together, these statements prove that $M_{T,h^{d}\cdot xy^{-a_1} \cdots xy^{-a_n}}$ is an $L$-space if and only if $d\in \{-1,0,-1\}.$ \\
 
\noindent \textbf{II.} Let us now consider the diffeomorphisms $\phi = h^d \cdot y^m$, where $m \in \mathbb{Z}$. We only need to deal with odd values of $d$ (when $d$ is even, $b_1(\ms)\geq1$, and, hence, $M_{T,\phi}$ is not an $L$-space). By Theorem \ref{thm:Tight}, $c(T,h \cdot y^m) \neq 0$. It suffices to show that $M_{T,h \cdot y^m}$ is an $L$-space since this would imply that $M_{T,h^{2d+1} \cdot y^m}$ is an $L$-space if and only if $d \in \{0,-1\}$ (which is the statement of Theorem \ref{thm:GOFClass}), by Corollary \ref{cor:Lspace}. Let us denote this 3-manifold by $Y(m)$. As Josh Greene points out, the manifolds $Y(m)$ are all rational surface singularities, and, hence, $L$-spaces by a theorem of N{\'e}methi in \cite{nem}. We give a self-contained proof below. The following lemma is very helpful.

\begin{lemma}
\label{lem:YRec}
$Y(m)$ an $L$-space $\Longleftrightarrow$ $Y(-m)$ an $L$-space $\Longleftrightarrow$ $Y(m+4)$ an $L$-space.
\end{lemma}

\begin{proof}[Proof of Lemma \ref{lem:YRec}]
 If $Y(m) = M_{T,h\cdot y^m}$ is an $L$-space, then Corollary \ref{cor:Lspace} implies that $M_{T,h^{-1}\cdot y^m}$ is an $L$-space, as above. This, in turn, implies that $Y(-m) = M_{T,h \cdot y^{-m}} \cong -M_{T,h^{-1} \cdot y^m}$ is an $L$-space, establishing the first equivalence. For the second equivalence, observe that closed braid specified by $(\x\y)^3\cdot \y^{m}$ is obtained by reversing the orientation on one of the components of the mirror of the closed braid specified by $(\x\y)^3\cdot \y^{-m-4}$. Hence, $M_{T,h\cdot y^{m}} \cong -M_{T,h\cdot y^{-m-4}}$; that is, $Y(m) \cong -Y(-m-4)$. Combined with the first equivalence, this shows that $Y(m)$ is an $L$-space if and only if $Y(m+4)$ is an $L$-space.
 \end{proof} 

Therefore, if $Y(0)$, $Y(-1)$, and $Y(-2)$ are $L$-spaces then $Y(m)$ is an $L$-space for all $m \in \mathbb{Z}$, completing the proof of Theorem \ref{thm:GOFClass}.(2). Recall that $Y(-2) = M_{T,h \cdot y^{-2}}$ is the double cover of $S^3$ branched along the closed 3-braid specified by $$(\x\y)^3 \cdot \y^{-2} = \x\y\x\y\x\y\cdot \y^{-2} = \x\y\y\x\y\y\cdot \y^{-2} = \x\y\y\x.$$ This braid has four crossings. It is therefore an alternating link, and, hence, $Y(-2)$ is an $L$-space. Likewise, $Y(-1) = M_{T,h \cdot y^{-1}}$ is the double cover of $S^3$ branched along the closed 3-braid specified by $$(\x\y)^3 \cdot y^{-1} = \x\y\x\y\x.$$ Since this link has five crossings, it too is alternating, and, hence, $Y(-1)$ is an $L$-space. This trick does not work for $Y(0)$; instead, we turn to a Heegaard diagram.

Below, we use the Heegaard diagram for $\ms$ constructed by Honda, Kazez, and Mati{'c} in \cite{hkm3}. Their construction is illustrated in Figure \ref{fig:Heeg}, which is a Heegaard diagram for $Y(0) = M_{T,h}$. Each of $T_{1/2}$ and $T_0$ is a copy of $T$. The Heegaard surface $\Sigma$ is the union of $T_{1/2}$ with $-T_0$ via an orientation reversing diffeomorphism of their boundaries. There are eight points of intersection between $\mathbb{T}_{\alpha}$ and $\mathbb{T}_{\beta}$ in $\text{Sym}^2(\Sigma)$. These are grouped according to their associated $\text{Spin}^c$ structures in the table below. To be clear, we have chosen an identification of $\text{Spin}^c(Y(0))$ with $H_1(Y(0);\mathbb{Z}) \cong \zzt \oplus \zzt.$

\begin{figure}[!htbp]
\begin{center}
\includegraphics[width=13cm]{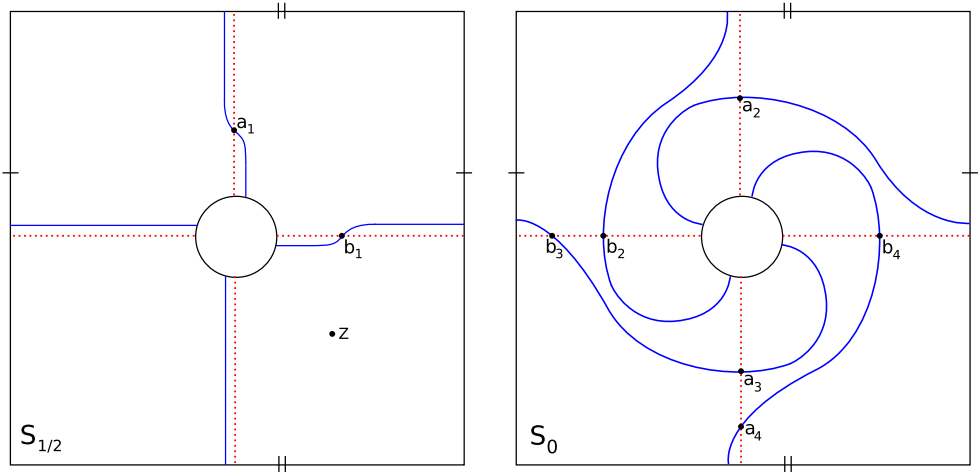}
\caption{A genus 2 Heegaard diagram for $Y(0)$. The $\alpha$ curves are dotted and the $\beta$ curves are solid. }
\label{fig:Heeg}
\end{center}
\end{figure}

\begin{table}[ht]
 \label{int}
\begin{center}
\begin{tabular}{|c|c|c|c|} \hline 
$(0,0)$ & $(1,0)$& $(0,1)$&$(1,1)$ \cr \hline \hline 
 $a_1 \times b_1$  & $a_1 \times b_3$ & $a_4 \times b_1$ &  $a_4 \times b_3$  \cr \hline 
 & $a_2 \times b_4$ & $a_2 \times b_2$ & \cr \hline 
 & $a_3 \times b_2$ & $a_3 \times b_4$ &  \cr \hline
\end{tabular}
\end{center}
\vspace{3mm}
 \caption{Generators of $\cf(Y(0))$ according to $\Sc$ structure.}
\end{table}

We wish to show that $Y(0)$ is an $L$-space; that is, $\hf(Y(0), \spc) \cong \zz$ for each $\spc\in \text{Spin}^c(Y(0)).$ This is evident when $\spc = (0,0)$ and $\spc =(1,1)$. Moreover, it is clear from the Heegaard diagram that in the chain complex $\cf(Y(0),(1,0))$, $a_3 \times b_2$ appears with multiplicity $\pm 1$ in $\partial (a_1\times b_3)$. Therefore, $\hf(Y(0),(1,0);\zt) \cong \zt$ for any field $\zt$, and hence, $\hf(Y(0),(1,0)) \cong \zz$. Likewise, in the complex $\cf(Y(0),(0,1))$, $a_3 \times b_4$ appears with multiplicity $\pm 1$ in $\partial (a_4\times b_1)$. Consequently, $\hf(Y(0),(0,1)) \cong \zz$. Thus, $Y(0)$ is an $L$-space. \\

%An alternate proof of the \emph{if} portion of this part of Theorem \ref{thm:GOFClass} can be given using the fact that $Y(m)$ can be written as the boundary of a 4-manifold corresponding to a plumbing diagram with only 1 bad vertex \cite{osz16}.

\noindent \textbf{III.} Finally, let us consider the diffeomorphisms $\phi = h^d \cdot x^m y^{-1},$ where $m \in \{-1,-2,-3\}$. Certainly $M_{T, x^{-m}y}$ is an $L$-space -- it is the branched double cover of $S^3$ along the braid given by $\x^{-m}\y$, which is the alternating link $-T(2,m)$. In fact, this branched cover is simply the lens space $L(m,1)$. Moreover, $c(T,x^{-m}y) \neq 0$ since any open book whose monodromy is the product of right-handed Dehn twists is Stein-fillable. Therefore, Corollary \ref{cor:Lspace} implies that $M_{T, h^{-2d}\cdot x^{-m}y}$ is an $L$-space if and only if $d \in \{0,1\}$. It follows that $M_{T, h^{2d}\cdot x^{m}y^{-1}}\cong - M_{T, h^{-2d}\cdot x^{-m}y}$ is an $L$-space if and only if $d\in \{0,1\}$.  

$M_{T, h\cdot x^my^{-1}}$ is also an $L$-space (this was the base case in the induction for the diffeomorphisms in Theorem \ref{thm:GOFClass}.(1)). Recall that $h\cdot x^{m}y^{-1} = x^{m+4}y$. When $m \in \{-1,-2,-3\}$, $x^{4+m}y$ is a composition of right-handed Dehn twists. Hence, $c(T,h\cdot x^my^{-1}) \neq 0$. As a result, Corollary \ref{cor:Lspace} implies that $M_{T,h^{2d+1}\cdot x^my^{-1}}$ is an $L$-space if and only if $d \in\{0,-1\}$. Thus, we have shown that $M_{T,h^{d} \cdot x^my^{-1}}$ is an $L$-space if and only if $d\in\{-1,0,1,2\}$.  

\end{proof}

It follows from Theorem \ref{thm:GOFClass} that any $\phi$ in Corollary \ref{cor:GOFClass} with $b_1(\ms)=0$ can be expressed as $\phi = h^{2d}\cdot\phi'$, where $M_{T,\phi'}$ is an $L$-space. That is to say, every 3-manifold of the form $M_{T,\phi}$ with $b_1(M_{T,\phi}) =0$ is obtained via $-1/d$-surgery on a genus one fibered knot in an $L$-space, as alluded to in Section \ref{section:tb.surg}. The corresponding torus bundles $\mt$ are obtained via $0$-surgeries on such knots.

\newpage
\section{A formula for the correction terms}
\label{section:corr}

In Section \ref{section:tb.msmt}, we compute the $\mathbb{Q}$-graded Heegaard Floer homologies $\hfp(M_{T,\phi},\so)$  and $\hfp(\mt,(\so,0))$ for all diffeomorphisms $\phi$ in Corollary \ref{cor:GOFClass} for which $b_1(M_{T,\phi})=0$. The only missing ingredients for an explicit description of these Floer homologies are, first, a precise formulation for the Floer homology of $1/n$- and $0$-surgeries on $\tref$, $-\tref$, and $\fige$, and, second, a formula for the correction terms $d(M_{T,\phi},\so)$. The first is provided by Ozsv{\'a}th and Szab{\'o} in \cite{osz6}, and the second consumes the rest of this section.

Before proceeding further, note that the Floer homology of $\ms$ is most interesting in the $\Sc$ structure $\so$. Indeed, as noted in the previous section, $\phi = h^{2d}\cdot\phi'$ for some $\phi'$ for which $M_{T,\phi'}$ is an $L$-space. It then follows from Proposition \ref{prop:Surgery} that for $\spc \neq \so$, $$\hfp(M_{T,\phi}, \spc) \cong \hfp(S^3_{-1/d}(\uknot))\{d(M_{T,\phi'}, \spc)\} \cong \tp_{d(M_{T,\phi'}, \spc)}.$$ If $\phi'$ is a diffeomorphism as in Corollary \ref{cor:GOFClass}.(1) or \ref{cor:GOFClass}.(3) it is possible to compute the correction term $d(M_{T,\phi'}, \spc)$ following \cite{osz12}, since, as we shall see Section \ref{section:tb.3braid}, $M_{T,\phi'}$ in this case is the branched double cover of $S^3$ along a quasi-alternating link (with a few exceptions for which the correction terms are already known). When $\phi'$ is a diffeomorphism as in Corollary \ref{cor:GOFClass}.(2), $M_{T,\phi'}$ is the boundary of a negative-definite plumbing with one bad vertex, so $d(M_{T,\phi'}, \spc)$ may be computed using the methods in \cite{osz16}. We use a different approach to compute the correction term in the $\text{Spin}^c$ structure $\so$.

\begin{proposition} 
\label{prop:Correction}
Below are the correction terms for several of the $L$-spaces among the manifolds $\ms$.
\begin{enumerate}

\item If $\phi = h \cdot xy^{-a_1} \cdots xy^{-a_n}$, where $a_i \geq 0$ and some $a_j \neq 0$ then $$d(\ms,\so) = (n+4 - \sum_{i=1}^n {a_i})/4.$$

\item If $\phi = xy^{-a_1} \cdots xy^{-a_n}$, where the $a_i \geq 0$, some $a_j \neq 0$, then $$d(\ms,\so)= (n - \sum_{i=1}^n {a_i})/4.$$ 

\item If $\phi = h\cdot y^m$, then $d(\ms,\so)= (m+4)/4.$

\item If $\phi = h\cdot x^my^{-1}$, for $m \in \{-1,-2,-3\},$ then $d(\ms,\so) = (m+3)/4.$

\item If $\phi = x^my^{-1}$, for $m \in \{-1,-2,-3\},$ then $d(\ms,\so) = (m+1)/4.$

\end{enumerate}

\end{proposition}

A key tool in the proof of Proposition \ref{prop:Correction} is the following result from \cite{osz1}.

\begin{proposition}
\label{prop:contact}
If $(S, \phi)$ is an open book decomposition for Y with binding $K$, and $\gamma \subset Y - K$ is a curve supported in a page of the open book decomposition, which is not homotopic within the page to the boundary, then $(S,  \gamma^{-1}\cdot \phi )$ induces an open book decomposition of $Y_{+1}(\gamma)$, and under the map

$$\widehat{F}_{W}: \hf(-Y) \longrightarrow \hf(-Y_{+1}(\gamma))$$

\noindent obtained by a 2-handle addition (and summing over all $spin^{c}$ structures on $W$), we have that $\widehat{F}_{W}(c(S, \phi)) = c(S, \gamma^{-1} \cdot \phi).$
\end{proposition}

Proposition \ref{prop:contact} implies the following.

\begin{lemma}
\label{lem:shift}
If $M_{T,\phi}$ and $M_{T, \gamma^{-1} \cdot \phi}$ are $L$-spaces, and $c(T, \gamma^{-1}\cdot \phi) \neq 0$ then $$d(M_{T, \gamma^{-1} \cdot \phi},\so) =  d(M_{T, \phi},\so) -1/4.$$
\end{lemma}

\begin{proof}[Proof of Lemma \ref{lem:shift}]
Recall that when $Y$ is an $L$-space, $d(Y, \spc)$ is simply the grading of the generator of $\hf(Y, \spc)$. Let $W$ be the cobordism from $-M_{T, \phi}$ to $-M_{T, \gamma^{-1} \cdot \phi}$ obtained from a 2-handle addition. Since $c(T, \gamma^{-1}\cdot \phi) \neq 0$, Proposition \ref{prop:contact} implies that $c(T,\phi) \neq 0$. Therefore, $c(T,\phi)$ generates $\hf(-M_{T, \phi},\so)$, while $c(T, \gamma^{-1}\cdot \phi)$ generates $\hf(-M_{T, \gamma^{-1} \cdot \phi},\so)$. Ozsv{\'a}th and Szab{\'o} show in \cite{osz5} that the diagram below commutes, where $i_1$ and $i_2$ are the maps on homology induced by the natural inclusion maps between chain complexes.

\begin{displaymath}
\xymatrix{
  \hfp(-M_{T, \phi},\so)  \ar[r]^{F^+_W}& \hfp(-M_{T, \gamma^{-1} \cdot \phi},\so)\\
 \hf(-M_{T, \phi},\so)  \ar[u]^{i_1}\ar[r]^{\widehat{F}_W}         &\hf(-M_{T, \gamma^{-1} \cdot \phi},\so)  \ar[u]^{i_2}  }
\end{displaymath} Since $-M_{T, \phi}$ and $-M_{T, \gamma^{-1} \cdot \phi}$ are $L$-spaces, the maps $i_1$ and $i_2$ are inclusions. Therefore, $F^{+}_W$ must be non-zero. Similarly, the maps $\pi_1$ and $\pi_2$ in the commutative diagram below, which are induced by the natural quotient maps between chain complexes, are surjective. Hence, $F^{\infty}_W$ is non-zero.

\begin{displaymath}
\xymatrix{
  \hfi(-M_{T, \phi},\so)  \ar[d]^{\pi_1}\ar[r]^{F^{\infty}_W}& \hfi(-M_{T, \gamma^{-1} \cdot \phi},\so) \ar[d]^{\pi_2} \\
 \hfp(-M_{T, \phi},\so)  \ar[r]^{F^+_W}         &\hfp(-M_{T, \gamma^{-1} \cdot \phi},\so)  }
\end{displaymath} Since $F^{\infty}_W$ is not identically zero it must be the case that $b_2^+(W) = 0$ \cite{osz5}. In fact, $W$ is negative-definite as both $-M_{T, \phi}$ and $-M_{T, \gamma^{-1} \cdot \phi}$ are rational homology 3-spheres.

For any $\text{Spin}^c$ structure $\spt$ on the cobordism $W$, $$F_{W,\bar{\spt}} = \mathfrak{J}F_{W,\spt} \mathfrak{J},$$ where $\mathfrak{J}:\hf(Y,\spc) \rightarrow \hf(Y,\bar{\spc})$ is induced by the isomorphism of chain complexes obtained by reversing the orientation of the Heegaard surface $\Sigma$ and switching the roles of the $\alpha$ and $\beta$ curves in a pointed Heegaard diagram $(\Sigma,\alpha,\beta,z)$ for the chain complex $\cf(Y)$ \cite{osz5}. Hence, if $\spt$ is a $\text{Spin}^c$ structure on $W$ which restricts to the $\text{Spin}^c$ structures $\so$ on $-M_{T, \phi}$ and $-M_{T, \gamma^{-1} \cdot \phi}$ then $$\widehat{F}_{W,\bar{\spt}}(c(T,\phi)) = \mathfrak{J}\widehat{F}_{W,\spt} \mathfrak{J}(c(T, \phi)) = \widehat{F}_{W,\spt} (c(T,\phi))$$ since the $\text{Spin}^c$ structures $\so$ are self-conjugate. 

If $\spt \neq \bar{\spt}$ then the contributions from the maps $\widehat{F}_{W,\spt}$ and $\widehat{F}_{W,\bar{\spt}}$ maps cancel (at least when working with $\zzt$ coefficients) when we sum over $\text{Spin}^c$ structures on $W$. The only contributions that survive are those which come from the maps $\widehat{F}_{W,\spt}$ where $\spt= \bar{\spt}$. For such $\spt $, $$c_1(\spt) = c_1(\bar{\spt})=-c_1(\spt) = 0.$$ Recall that the grading shift of the map $\widehat{F}_{W, \spt}$ is $$\frac{c_1(\spt)^2 -2\chi(W) -3\sigma(W)}{4}.$$ When $\spt =\bar{\spt}$ this grading shift is $1/4$ since $c_1(\spt)^2=0$, $\chi(W) = 1$, and $\sigma(W) = -1$. Therefore, $$d(-M_{T, \gamma^{-1} \cdot \phi},\so) =  d(-M_{T, \phi},\so) +1/4.$$ In general, $d(-Y, \spc) = -d(Y, \spc)$. Hence, $$d(M_{T, \gamma^{-1} \cdot \phi},\so) =  d(M_{T, \phi},\so) -1/4.$$
\end{proof}

\begin{proof}[Proof of Proposition \ref{prop:Correction}]
For $\phi$ in Proposition \ref{prop:Correction}.(1) or \ref{prop:Correction}.(3) - \ref{prop:Correction}.(5), the proof follows directly from Lemma \ref{lem:shift}. A bit more work is required for $\phi$ in Proposition \ref{prop:Correction}.(2). \\

\noindent \textbf{I.} Let $\phi = h \cdot xy^{-a_1} \cdots xy^{-a_n}$, where the $a_i \geq 0$ and some $a_j \neq 0$. $M_{T,y\cdot\phi} \cong M_{T, h \cdot xy^{-a_1} \cdots xy^{-a_n+1}}$ and $M_{T,\phi} \cong M_{T, h \cdot xy^{-a_1} \cdots xy^{-a_n}}$ are both $L$-spaces, and $c(T, \phi) \neq 0$ by Theorem \ref{thm:GOFClass} and Theorem \ref{thm:Tight}. Lemma \ref{lem:shift} implies that $$d(M_{T, h \cdot xy^{-a_1} \cdots xy^{-a_n}} ,\so) = d(M_{T, h \cdot xy^{-a_1} \cdots xy^{-a_n+1}}, \so) -1/4.$$ Iterating this, we find that \begin{equation} \label{eqn:1} d(M_{T,h \cdot xy^{-a_1} \cdots xy^{-a_n}},\so) = d(M_{T,h\cdot x^n y^{-1}}, \so) - (\sum_{i=1}^n{a_i}-1)/4.\end{equation} From the previous section, $M_{T, h\cdot x^n y^{-1}} \cong L(n+4,1)$. For $n$ odd, Ozsv{\'a}th and Szab{\'o} show that $L(n+4,1)$ has a unique self-conjugate $\text{Spin}^c$ structure $\so$ for which $d(L(n+4,1),\so)=(n+3)/4$ \cite{osz6}. By Lemma \ref{lem:shift}, $$d(M_{T, h\cdot x^ny^{-1}} ,\so) = (n+3)/4$$ for any $n>0$. Plugging this into Equation \ref{eqn:1}, we see that $$d(M_{T,h \cdot xy^{-a_1} \cdots xy^{-a_n}},\so) = (n+3)/4 -(\sum_{i=1}^n{a_i}-1)/4 = (n+4 -\sum_{i=1}^n{a_i})/4.$$ \\

\noindent \textbf{II.} Let $\phi = xy^{-a_1} \cdots xy^{-a_n}$, where the $a_i \geq 0$ and some $a_j \neq 0$. We abuse notation in the following way. If $K$ is a knot in $Y$, and $\gamma \subset Y$ is a knot which does not link $K$ then we also denote by $K$ the corresponding knot in $Y_{p/q}(\gamma)$. Now, the manifolds $$-M_{T,xy^{-a_1} \cdots xy^{-a_n+1}}, \,-M_{T,xy^{-a_1} \cdots xy^{-a_n}},\, \text{and } Q(a_1+1,a_2,\dots,a_{n-1},a_{n-1}+1)$$ are obtained via $\infty$-, $-1$-, and $0$-surgeries on a knot isotopic to a copy of the curve $y$ on a page of the open book decomposition $(T,xy^{-a_1} \cdots xy^{-a_n+1})$. As a result, there is a surgery exact triangle on knot Floer homology,
\begin{displaymath}
\xymatrix{
  \hfk(-M_{T,xy^{-a_1} \cdots xy^{-a_n+1}},K)  \ar[r]^{\widehat{F}_W}&  \hfk(-M_{T,xy^{-a_1} \cdots xy^{-a_n}},K) \ar[d] \\
 &\hfk(Q(a_1+1,a_2,\dots,a_{n-1},a_{n-1}+1),K) \ar[ul] }
\end{displaymath}
where $K$ is the binding of the open book $(T,xy^{-a_1} \cdots xy^{-a_n+1})$. 

Since the curve $y$ does not link $K$, these maps preserve the Alexander filtrations. Moreover, it is easy to see from the surgery diagram for $Q(a_1+1,a_2,\dots,a_{n-1},a_{n-1}+1)$ that the induced knot $K \subset Q(a_1+1,a_2,\dots,a_{n-1},a_{n-1}+1)$ has genus 0. Therefore, the map $$\widehat{F}_W:\hfk(-M_{T,xy^{-a_1} \cdots xy^{-a_n+1}},K,-1) \rightarrow \hfk(-M_{T,xy^{-a_1} \cdots xy^{-a_n}},K,-1)$$ is an isomorphism. Iterating, we obtain an isomorphism $$\widehat{F}_{\ww}: \hfk(-M_{T,x^ny^{-1}},K,-1) \rightarrow \hfk(-M_{T,xy^{-a_1} \cdots xy^{-a_n}},K,-1),$$ where $\widehat{F}_{\ww}$ is the composition of $\sum_{i=1}^n{a_i}-1$ maps like $\widehat{F}_W$.

Likewise, there is a surgery exact triangle relating the knot Floer homologies $$\hfk(-M_{T,xy^{-a_1} \cdots xy^{-a_n}},K),\text{ and }\hfk(-M_{T,y^{-a_1} \cdots xy^{-a_n}},K),$$ where the second 3-manifold is obtained from the first via $-1$-surgery around a knot isotopic to a copy of $x$ on a page of the open book decomposition $(T,xy^{-a_1} \cdots xy^{-a_n})$. As before, the induced map $$\widehat{F}_{W'}:\hfk(-M_{T,xy^{-a_1} \cdots xy^{-a_n}},K,-1) \rightarrow \hfk(-M_{T,y^{-a_1} \cdots xy^{-a_n}},K,-1)$$ is an isomorphism. We iterate again to obtain an isomorphism $$\widehat{F}_{\ww'}: \hfk(-M_{T,xy^{-a_1} \cdots xy^{-a_n}},K,-1) \rightarrow \hfk(-M_{T,xy^{-\sum{a_i}}},K,-1),$$ where $\widehat{F}_{\ww'}$ is the composition of $n-1$ maps like $\widehat{F}_{W'}$. 

The composition $$\widehat{F}_{\ww'}\circ \widehat{F}_{\ww}: \hfk(-M_{T,x^ny^{-1}},K,-1) \rightarrow \hfk(-M_{T,xy^{-\sum{a_i}}},K,-1)$$ is therefore an isomorphism. In Subsection 5.2 of \cite{jabmark}, Jabuka and Mark show that the knot Floer homology $\hfk(-M_{T,x^ny^{-1}},K,-1)$ is supported in grading $(-n-3)/4$, and $\hfk(-M_{T,xy^{-\sum{a_i}}},K,-1)$ is supported in grading $(\sum{a_i} -5)/4$. Thus, $\widehat{F}_{\ww'}\circ \widehat{F}_{\ww}$ is a map of degree $(\sum_{i=1}^n{a_i}+n-2)/4.$ Each of the maps $\widehat{F}_W$ which goes into the composition $\widehat{F}_{\ww}$ is a sum $$\widehat{F}_W=\sum_{\spt\in \text{Spin}^c(W)}{\widehat{F}_{W, \spt}}.$$ If $W$ is the cobordism from $-M_{T,xy^{-a_1} \cdots xy^{-a_n+1}}$ to $-M_{T,xy^{-a_1} \cdots xy^{-a_n}}$ then it is shown in \cite{bald1} that $$\widehat{F}_W: \hf(-M_{T,xy^{-a_1} \cdots xy^{-a_n+1}}) \rightarrow \hf(-M_{T,xy^{-a_1} \cdots xy^{-a_n}})$$ is injective. Therefore, $W$ is negative definite and, hence, the degree shift of $\widehat{F}_{W,t}$ is at most $1/4$ for every $\text{Spin}^c$ structure $\spt $ on $W$. Precisely the same analysis can be carried out for the maps $\widehat{F}_{W', \spt}$ (in this case, the maps $\widehat{F}_{W'}$ on $\hf$ are surjective which also implies that $W'$ is negative-definite). 

Since there are $\sum_{i=1}^n{a_i} -1$ maps like $\widehat{F}_W$ and $n-1$ maps like $\widehat{F}_{W'}$ going into the composition $\widehat{F}_{\mathbb{W}'} \circ \widehat{F}_{\mathbb{W}}$, we can conclude that the map $\widehat{F}_{\mathbb{W}'} \circ \widehat{F}_{\mathbb{W}}$ has degree at most $$(\sum_{i=1}^n{a_i} -1 +n-1)/4 = (\sum_{i=1}^n{a_i} +n-2)/4.$$ On the other hand, we have found that $\widehat{F}_{\mathbb{W}'} \circ \widehat{F}_{\mathbb{W}}$ has degree \emph{exactly} $(\sum_{i=1}^n{a_i} +n-2)/4$. Therefore, each of the maps $\widehat{F}_W$ and $\widehat{F}_{W'}$ in the composition $\widehat{F}_{\mathbb{W}'} \circ \widehat{F}_{\mathbb{W}}$ has degree $1/4$. 

In each of the knot Floer homologies above, the generator in Alexander grading $-1$ is supported in the $\text{Spin}^c$ structure $\so$. It follows that there is a $\text{Spin}^c$ structure $\spt $ on the cobordism $\ww' \circ \ww$ restricting to $\so$ on $-M_{T,x^ny^{-1}}$ and on $-M_{T,xy^{-\sum{a_i}}}$ for which $c_1(\spt)^2 =0$. Thus, $$F^{\infty}_{\ww' \circ \ww,t}: \hfi(-M_{T,x^ny^{-1}},\so) \rightarrow \hfi(-M_{T,xy^{-\sum{a_i}}},\so)$$ is an isomorphism of degree $(\sum_{i=1}^n{a_i} +n-2)/4$, and the map $F^{+}_{\ww' \circ \ww,t}$ induced on $\hfp$ is surjective. However, since $$d(-M_{T,x^ny^{-1}},\so) -d(-M_{T,xy^{-\sum{a_i}}},\so) = (\sum_{i=1}^n{a_i} +n-2)/4$$ it follows that the map $F^{+}_{\ww' \circ \ww, \spt}$ and, hence, the map $$\widehat{F}_{\ww' \circ \ww, \spt}: \hf(-M_{T,x^ny^{-1}},\so) \rightarrow \hf(-M_{T,xy^{-\sum{a_i}}},\so)$$ is an isomorphism (we are using the fact that our manifolds are $L$-spaces). Therefore, each of the maps $$\widehat{F}_{W, \spt |_W}: \hf(-M_{T,xy^{-a_1} \cdots xy^{-a_n+1}},\so) \rightarrow \hf(-M_{T,xy^{-a_1} \cdots xy^{-a_n}},\so)$$ is an isomorphism of degree $1/4$, and by induction, $$d(-M_{T,xy^{-a_1} \cdots xy^{-a_n}},\so) = d(-M_{T,x^ny^{-1}},\so) + (\sum_{i=1}^n{a_i}-1)/4 = (-n + \sum_{i=1}^n{a_i})/4.$$ Finally, $$d(M_{T,xy^{-a_1} \cdots xy^{-a_n}},\so) = -d(-M_{T,xy^{-a_1} \cdots xy^{-a_n}},\so) = (n - \sum_{i=1}^n{a_i})/4.$$ 

In the analysis above, we stopped just short of computing the Maslov grading of the generator of $\hfk(M_{T,xy^{-a_1}\cdots xy^{-a_n}},K,-1,\so)$, but it is easy to see that this grading is one less than the Maslov grading of the generator of $\hf(M_{T,xy^{-a_1}\cdots xy^{-a_n}},\so)$. This validates the comments following Remark \ref{remark:ss} regarding the relative Maslov gradings on $\hfk(\ms,K,\so)$ for $c(T,\phi)=c(T,\phi^{-1})=0$.\\

\noindent \textbf{III.} Let $\phi = h\cdot y^m$. These $\ms$ are all $L$-spaces, and $c(T,\phi) \neq 0$. Hence, Lemma \ref{lem:shift} implies that $$d(M_{T,h\cdot y^{m-1}},\so)=d(M_{T,h\cdot y^{m}},\so) -1/4.$$ When $m=-1,$ $$M_{T,\phi} \cong M_{T,xyxyxy\cdot y^{-1}} \cong M_{T,xxyxx} \cong L(4,1).$$ Therefore, $d(M_{T,h\cdot y^{-1}})=3/4$, and it follows by induction that $d(M_{T,h \cdot y^m},\so) = (m+4)/4.$ \\

\noindent \textbf{IV.} Let $\phi = h\cdot x^my^{-1}$, for $m \in \{-1,-2,-3\}.$ As we have seen, $M_{T,h\cdot x^my^{-1}} \cong L(m+4,1)$. Therefore, $d(M_{T,\phi},\so) = (m+3)/4.$\\

\noindent \noindent \textbf{V.} Let $\phi = x^my^{-1}$, for $m \in \{-1,-2,-3\}$. Then $M_{T,\phi} \cong -L(-m,1) \cong L(m,1)$. Hence, $d(M_{T,\phi},\so) = -(-m-1)/4 = (m+1)/4.$

\end{proof}
\newpage

\section{The $\mathbb{Q}$-graded Floer homology of $\ms$ and $\mt$}
\label{section:tb.msmt}
\subsection{$\ms$, $\mt$, and embedded contact homology}
In \cite{osz6}, Ozsv{\'a}th and Szab{\'o} compute the Heegaard Floer homology of $1/n$-surgery on $\tref$ for all $n \in \zz$, and of $1/n$-surgery on $\fige$ for $n >0$. We have copied and expanded their results using the symmetry of Heegaard Floer homology under orientation reversal (see \cite{osz5}) in the following proposition. Recall that the subscripts indicate grading.

\begin{proposition}
\label{prop:TE} Below are the Heegaard Floer homologies of $1/n$-surgeries on $\tref$, $-\tref$, and $\fige$.
\begin{itemize}
\item $\hfp(S^3_{1/n}(\tref)) \cong 
\left\{\begin{array}{lr} 
\tp_{-2}\oplus \zz^{n-1}_{-2}, &{\text{if $n>0$}}\\ 
\tp_{0}\oplus \zz^{-n}_{-1}, &{\text{if $n\leq0$}}
\end{array}\right.$\\

\item $\hfp(S^3_{1/n}(-\tref)) \cong 
\left\{\begin{array}{lr} 
\tp_{0}\oplus \zz^{n}_{0}, &{\text{if $n\geq0$}}\\ 
\tp_{2}\oplus \zz^{-n-1}_{1}, &{\text{if $n<0$}}
\end{array}\right.$\\

\item $\hfp(S^3_{1/n}(\fige)) \cong 
\left\{\begin{array}{lr} 
\tp_{0}\oplus \zz^{n}_{-1}, &{\text{if $n\geq0$}}\\ 
\tp_{0}\oplus \zz^{-n}_{0}, &{\text{if $n<0$}}
\end{array}\right.$\\
\end{itemize}
\end{proposition}

Combined with Proposition \ref{prop:TE}, the following theorem gives an explicit description of the Heegaard Floer homology of any 3-manifold $\ms$ with $b_1(\ms)=0$.

\begin{theorem} Suppose that $b_1(\ms)=0$. For $\spc =\so$,
\label{thm:Grading}
\begin{enumerate}
\item If $\phi = h^d \cdot xy^{-a_1} \cdots xy^{-a_n}$, where the $a_i \geq 0$ and some $a_j \neq 0$, then $$\hfp(M_{T,\phi},\so) \cong
\left\{\begin{array}{ll} 
\hfp(S^3_{-1/k}(\tref))\{(n+4 - \sum_{i=1}^n {a_i})/4\}, &{\text{if $d=2k+1$}}\\
\hfp(S^3_{-1/k}(\fige))\{(n - \sum_{i=1}^n {a_i})/4\}, &{\text{if $d=2k$.}}
\end{array}\right.$$ 

\item If $\phi = h^{2d+1} \cdot y^m$, for $m \in \mathbb{Z}$, then
$$\hfp(M_{T,\phi},\so) \cong \hfp(S^3_{-1/d}(\tref))\{ (m+4)/4\}.$$ 

\item If $\phi = h^{d} \cdot x^m y^{-1},$ where $m \in \{-1,-2,-3\}$, then
$$\hfp(M_{T,\phi},\so) \cong
\left\{\begin{array}{ll} 
\hfp(S^3_{-1/k}(\tref))\{(m+3)/4\}, &{\text{if $d=2k+1$}}\\
\hfp(S^3_{-1/k}(-\tref))\{(m+1)/4\}, &{\text{if $d=2k$.}}
\end{array}\right.$$ 
\end{enumerate}
For all $\spc \neq \so$, $$\hfp(\ms, \spc) \cong \tp_{d(\ms, \spc)}.$$
\end{theorem}

\begin{proof}[Proof of Theorem \ref{thm:Grading}]
This follows directly from Propositions \ref{prop:Surgery} and \ref{prop:Correction}, and Theorems \ref{thm:GOFClass} and \ref{thm:Tight}. 
\end{proof}

As mentioned in Section \ref{section:tb.surg}, $\text{Spin}^c(\mt)$ can be identified with $\text{Spin}^c(\ms) \oplus \zz$, where the first Chern class of the $\text{Spin}^c$ structure $(\spc,i)$ evaluates $2i$ on the fiber $\widehat{T}$. Since $g(\widehat{T})=1$, the adjunction inequality in Heegaard Floer homology implies that $\hfp(\mt)$ is entirely supported in $\text{Spin}^c$ structures of the form $(\spc,0)$ \cite{osz14}. 

In \cite{osz6}, Ozsv{\'a}th and Szab{\'o} compute the Heegaard Floer homology of $0$-surgery on $\tref$ and $\fige$. We have reproduced their results below.

\begin{proposition}
\label{prop:TE0} Below are the Heegaard Floer homologies of $0$-surgery on $\tref$, $-\tref$, and $\fige$ in the unique $\text{Spin}^c$ structures in which they are supported.
\begin{itemize}
\item $\hfp(S^3_0(\tref),0) \cong \tp_{-\frac{1}{2}}\oplus \tp_{-\frac{3}{2}},$
\item $\hfp(S^3_0(-\tref),0) \cong \tp_{\frac{3}{2}}\oplus \tp_{\frac{1}{2}},$
\item $\hfp(S^3_0(\fige),0) \cong \tp_{\frac{1}{2}}\oplus \tp_{-\frac{1}{2}} \oplus \zz_{-\frac{1}{2}}.$
\end{itemize}
\end{proposition}

Combined with the following theorem, Proposition \ref{prop:TE0} gives a complete description of the Heegaard Floer homology of the torus bundles $\mt$ with $b_1(\mt)=1$.

\begin{theorem} Suppose that $b_1(\mt)=1$. For $\spc =\so$,
\label{thm:TB}
\begin{enumerate}
\item If $\phi = h^d \cdot xy^{-a_1} \cdots xy^{-a_n}$, where the $a_i \geq 0$ and some $a_j \neq 0$, then $$\hfp(\mt,(\so,0)) \cong
\left\{\begin{array}{ll} 
\hfp(S^3_{0}(\tref))\{(n+4 - \sum_{i=1}^n {a_i})/4\}, &{\text{if $d$ is odd}}\\
\hfp(S^3_{0}(\fige))\{(n - \sum_{i=1}^n {a_i})/4\}, &{\text{if $d$ is even.}}
\end{array}\right.$$ 

\item If $\phi = h^{2d+1} \cdot y^m$, for $m \in \mathbb{Z}$, then
$$\hfp(\mt,(\so,0)) \cong \hfp(S^3_{0}(\tref))\{ (m+4)/4\}.$$ 

\item If $\phi = h^{d} \cdot x^m y^{-1},$ where $m \in \{-1,-2,-3\}$, then
$$\hfp(\mt,(\so,0)) \cong
\left\{\begin{array}{ll} 
\hfp(S^3_{0}(\tref))\{(m+3)/4\}, &{\text{if $d$ is odd}}\\
\hfp(S^3_{0}(-\tref))\{(m+1)/4\}, &{\text{if $d$ is even.}}
\end{array}\right.$$ 
\end{enumerate}
For all $\spc\neq \so$, $$\hfp(\mt,(\spc,0))\cong \tp_{d(\ms, \spc)+1/2} \oplus \tp_{d(\ms, \spc)-1/2}.$$ For all $(\spc,i)$, where $i\neq 0$, $$\hfp(\mt,(\spc,i))=0.$$
\end{theorem}

\begin{proof}[Proof of Theorem \ref{thm:TB}] This proof is identical to that of Theorem \ref{thm:Grading}.
\end{proof}

Below, we compare Heegaard Floer homology and embedded contact homology for torus bundles with pseudo-Anosov monodromy. Embedded contact homology associates a graded Abelian group $ECH_*(Y,\lambda;\Gamma)$ to a closed oriented 3-manifold $Y$, a contact form $\lambda$, and a homology class $\Gamma \in H_1(Y)$ \cite{hs}. $ECH$ is conjecturally independent of the contact form $\lambda$. In fact, $ECH(Y,\lambda;\Gamma)$ is conjectured to be isomorphic to the Seiberg-Witten Floer homology, $HM^{\vee}_*(-Y,s_{\xi}(\Gamma)),$ and to $\hfp(-Y,s_{\xi}(\Gamma))$, where $$s_{\xi}: H_1(Y)\rightarrow \Sc(Y)$$ is an affine bijection which sends $0$ to the $\Sc$ structure associated to the contact 2-plane field $\xi=\text{ker}(\lambda)$.

Observe that $$H_1(\mt;\zz)\cong \frac{H_1(\widehat{T})}{\text{im}(id-\phi_*)} \oplus \zz,$$ where $\phi_*$ is the map induced on $H_1(\widehat{T})$ by $\phi$. Therefore, we may denote a homology class in $H_1(\mt)$ by $(\gamma,i)$. In his thesis, Eli Lebow constructs a sequence of contact forms $\{\lambda_n\}$ on $\mt$ for any pseudo-Anosov $\phi$ (recall that these are the diffeomorphisms $\phi$ in Corollary \ref{cor:GOFClass}.1), and he shows that $ECH(\mt,\lambda_n;\Gamma)$ is independent of $n$ \cite{lebow}. A quick comparison of Theorem 1.1 in Lebow's thesis with Theorem \ref{thm:TB} above shows that the following are isomorphic as relatively graded $\zz$-modules:

$$ECH(\mt,\lambda_n;(\gamma,i)) \cong
\left\{\begin{array}{ll} 
\hfp(-\mt,(\spc,i)), &{\text{for any $\spc $, if $i\neq 0$,}}\\
\hfp(-\mt,(\spc,0)), &{\text{for any $\spc\neq \so$, if $i=0$ and $\gamma \neq 0$,}}\\
\hfp(-\mt,(\so,0)), &{\text{if $i=0$ and $\gamma = 0$.}}
\end{array}\right.$$ 

It would be interesting to show that the fully twisted versions of both theories are isomorphic as relatively graded $(\zz[U]\otimes \wedge^*(H_1(Y)/\text{Tors}) )$-modules, as is done in \cite{hs} by Hutchings and Sullivan in the case that $\phi$ is the identify.

\subsection{Other surface bundles over $S^1$}
Let $\Sigma_g$ denote the closed surface of genus $g>1$. In \cite{jabmark,jabmark2}, Jabuka and Mark compute the Heegaard Floer homology in non-torsion $\Sc$ structures of $\Sigma_g$-bundles over $S^1$ with monodromy $\phi$, where $\phi$ is one of the following products of Dehn twists:\begin{itemize}
\item $\phi = \sigma^n$, where $\sigma$ is a genus 1 separating curve on $\Sigma_g$
\item $\phi = \alpha^n \beta^m$, where $\alpha$ and $\beta$ are dual, non-separating curves on $\Sigma_g$.
\end{itemize}

Using their notation througout this subsection, we denote such an $\Sigma_g$-bundle by $M(\phi)$. When $n$ and $m$ are both non-zero, Jabuka and Mark compute $\hfp(M(\alpha^n \beta^m),\spc)$ directly from the knot Floer homology $\hfk(M_{T,x^ny^m},K,\spt)$ and the associated $d^1$ differentials (where $K$ is the binding of the open book decomposition $(T,x^ny^m)$, and $\spt$ is any $\Sc$ structure on $M_{T,x^ny^m}$). We can therefore extend their results without any additional effort using our computations of knot Floer homology and $d^1$ differentials in Proposition \ref{prop:HFK}, combined with our classification of $L$-spaces among the manifolds $M_{T,\phi}$ in Theorem \ref{thm:GOFClass}. Let us recall their setup. 

\begin{figure}[!htbp]
\begin{center}
\includegraphics[width=10cm]{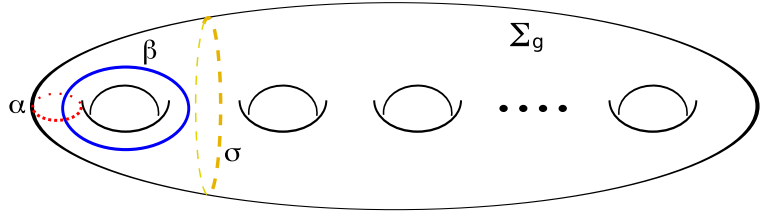}
\caption{The surface $\Sigma_g$.}
\label{fig:Sigmag}
\end{center}
\end{figure}

To begin with, observe that $H_2(M(\phi);\zz)\cong \zz \oplus \text{ker}(id-\phi_*).$ Let $\mathfrak{S}_k$ denote the set of $\spc\in\Sc(M(\phi))$ which satisfy the following conditions:
\begin{itemize}
\item $\langle c_1(\spc),[\Sigma_g] \rangle = 2k$, and
\item $\langle c_1(\spc),[\gamma] \rangle = 0$ for all $[\gamma]$ coming from $H_1(\Sigma_g)$.
\end{itemize}
These conditions specify the $\Sc$ structure $\spc$ up to torsion. Note that for $\spc \in \mathfrak{S}_k$, the adunction inequality implies that $\hfp(M(\phi), \spc)=0$ unless $|k|\leq g-1$. 

If $G$ is a graded group, let $G_{(m)}$ denote homogeneous piece of $G$ in grading $m$. Let $X(g,d)$ be the $\zz$-graded group defined by $X(g,d)_{(j)} = H_{g-j}(\text{Sym}^d(\Sigma_g);\zz)$ when $d\geq0$ (see \cite{imac} for an explicit description of the homology of a symmetric product). $X(g,d)$ is defined to be zero for $d<0$. For any diffeomorphism $\phi$ of $T$ written as a product of Dehn twists around $x$ and $y$, we let $\phi_{\Sigma}$ denote the diffeomorphism of $\Sigma_g$ obtained from $\phi$ by substituting $\alpha$ for $x$ and $\beta$ for $y$. We may now state a generalization of Theorem 1.2 from \cite{jabmark}.

\begin{theorem}
\label{thm:Jab}
Let $k$ be an integer with $0<|k|\leq g-1$, and define $d=g-1-|k|$. Let $\phi$ be a diffeomorphism of $T$ written as a product of Dehn twists in $x$ and $y$. If $\ms$ is an $L$-space, then there is a unique $\Sc$ structure $\spc_k\in \mathfrak{S}_k$ on the $\Sigma_g$-bundle $M(\phi_{\Sigma})$ for which there is an isomorphism of relatively graded groups

$$\hfp(M(\phi_{\Sigma}), \spc_k) \cong
\left\{\begin{array}{ll} 
X(g-1,d-1)\oplus \wedge^{2g-2-d}H^1(\Sigma_{g-1})_{(g-1-d)}\\
X(g-1,d-1)\{-1\}\oplus \wedge^{2g-2-d}H^1(\Sigma_{g-1})_{(g-1-d)}\\
X(g-1,d-1)\{-2\}\oplus \wedge^{2g-2-d}H^1(\Sigma_{g-1})_{(g-1-d)}.
\end{array}\right.$$ 
The three lines on the right-hand side above correspond to the three cases $c(T,\phi)\neq 0$, $c(T,\phi)=c(T,\phi)= 0$, and $c(T,\phi^{-1})\neq 0$, respectively. For all other $\spc\in \mathfrak{S}_k$, there is an isomorphism $$\hfp(M(\phi_{\Sigma}), \spc) \cong X(g-1,d-1).$$

\end{theorem}

\begin{proof}
Jabuka and Mark's Theorem 1.2 gives the Heegaard Floer homology of $\Sigma_g$-bundles of the form $M(\al^n\be^m)$ for $n$ and $m$ both non-zero. Theorem \ref{thm:Jab} is nearly a verbatim copy of Theorem 1.2, only we have replaced the conditions $m,n>0$, $m\cdot n <0$, and $m,n<0$ with $c(T,\phi)\neq 0$, $c(T,\phi)=c(T,\phi)= 0$, and $c(T,\phi^{-1})\neq 0$, respectively. As mentioned above, Jabuka and Mark's proof of Theorem 1.2 relies only on the knot Floer homology $\hfk(M_{T,x^ny^m},K,\spt)$ and the associated $d^1$ differentials in the spectral sequence for $K\subset M_{T,x^ny^m}$. Since the structure of these knot Floer homologies and differentials are preserved under the replacements we have made (see Proposition \ref{prop:HFK}), our more general statement follows directly from the proof of Theorem 1.2.
\end{proof}

\begin{figure}[!htbp]
\begin{center}
\includegraphics[width=10cm]{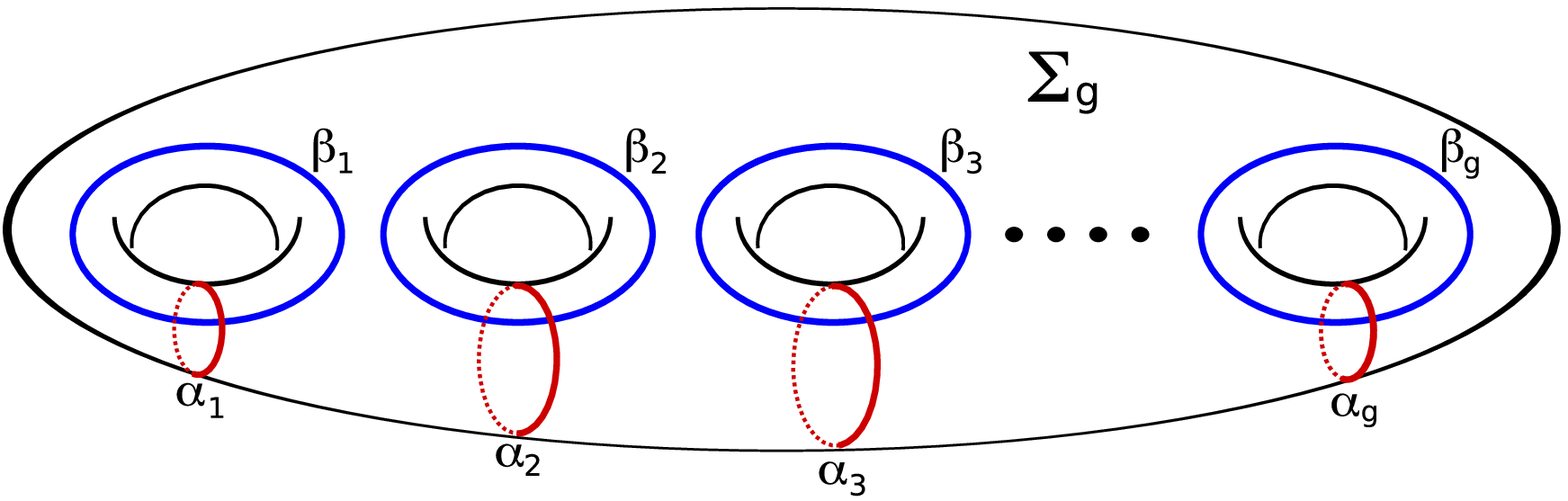}
\caption{}
\label{fig:Sigmag2}
\end{center}
\end{figure}

It should be possible, using the results in this paper, to compute the Heegaard Floer homology of more general surface bundles of the form $M(\phi)$, where $\phi$ is any composition of the Dehn twists around the curves $\alpha_1,\beta_1,\dots, \alpha_g,\beta_g$ depicted in Figure \ref{fig:Sigmag2}. We will return to this in a future paper.

There are other Floer theories associated to surface diffeomorphisms of $\Sigma_g$ (see \cite{eftek, hs2, seidel}) which appear to be linked with Heegaard Floer homology. For instance, it is conjectured that the Lagrangian Floer homology, $HF(\psi)$, of a diffeomorphism $\psi$ agrees with $\hfp(M(\psi), \spc_{k-2})$. In \cite{jabmark}, Jabuka and Mark verify this conjecture in all cases for which both theories have been computed by comparing their results with those of Eftekhary and Seidel \cite{eftek,seidel}. In addition, Jabuka and Mark show that $\hfp(M(\psi))$ agrees with the periodic Floer homology, $HP(\psi)$, of Hutchings and Sullivan (\cite{hs2}) where both have been calculated.

It would be interesting to use our extension of Jabuka and Mark's Theorem 1.2 to verify these conjectural relationships for more of the surface diffeomorphisms $\phi_{\Sigma}$. 
\newpage

\section{Stein-fillings of $\xi_{T,\phi}$}
\label{section:stein}

In \cite{bald1}, we prove that a diffeomorphism $\phi = x^{a_1}y^{b_1}\cdots x^{a_n}y^{b_n}$ is isotopic to the product of $k$ right-handed Dehn twists around non-separating curves in $T$ only if $$\sum_{i=1}^n a_i+b_i = k.$$ Along with Theorem \ref{thm:Tight}, this shows that any $\phi = h^d\cdot xy^{-a_1}\cdots xy^{-a_n}$ with $d>0$ and $6d+n-\sum a_i < 0$ is in $\Tight-\Dehn$. Honda, Kazez, and Mati{\'c} prove a very similar statement in \cite{hkm2}. 

Note that $\phi \in \Stein$ does not \emph{a priori} imply that $\phi \in \Dehn$. Indeed, Giroux proves that a Stein-fillable contact structures is compatible with \emph{some} open book whose monodromy is the product of right-handed Dehn twists \cite{giroux}, but this open book need not be of the form $(T,\phi)$. In this section, we produce an infinite family of diffeomorphisms $\phi \in \Tight -\Stein$ using the $\mathbb{Q}$-gradings on the Heegaard Floer homology groups $\hf(\ms,\so)$. Below is the main theorem of this section.

\begin{theorem}
\label{thm:stein}
If $\ms$ is an $L$-space and $(W,J)$ is a Stein-filling of $\xi_{T,\phi}$ then $$d(\ms,\so) = \frac{\chi(W)-1}{4} \geq 0.$$
\end{theorem}

For example, if $\phi = h \cdot xy^{-a_1}\cdots xy^{-a_n}$ then $d(\ms,\so) = (4+n-\sum a_i)/4$ by Proposition \ref{prop:Correction}. Therefore, $\xi_{T,\phi}$ is tight (by Theorem \ref{thm:Tight}), but not Stein-filllable as long as $(4+n-\sum a_i)/4<0$. 

Recall from \cite{osz1} that if $(Y,\xi)$ is a contact 3-manifold with Stein-filling $(W,J)$ then the grading of $c(\xi)\in \hf(-Y)$ is given by \begin{equation} \label{eqn:grading} \text{gr}(c(\xi)) = -\frac{(c_1(W,J))^2-2\chi(W)-3\sigma(W) +2}{4}.\end{equation} Furthermore, if $Y$ is an $L$-space then $c(\xi)$ generates $\hf(-Y,\spc_{\xi}),$ and Equation \ref{eqn:grading} implies that \begin{equation}\label{eqn:grading2} d(Y,\spc_{\xi})=-d(-Y,\spc_{\xi}) = \frac{(c_1(W,J))^2-2\chi(W)-3\sigma(W) +2}{4}.\end{equation}

\begin{proof}[Proof of Theorem \ref{thm:stein}]
Suppose that $(W,J)$ is a Stein-filling of the contact structure $\xi_{T,\phi}$ and that $\ms$ is an $L$-space. Suppose, for a contradiction, that $J \neq \bar{J}$. Then $(W,\bar{J})$ is a Stein-filling of $\overline{\xi_{T,\phi}}$, and Plamenevskaya proves in \cite{pla3} that $c(\xi_{T,\phi}) \neq c(\overline{\xi_{T,\phi}})$. Since the $\Sc$ structure $\so$ associated to the contact structure $\xi_{T,\phi}$ is self-conjugate, the contact invariants $c(\xi_{T,\phi})$ and $c(\overline{\xi_{T,\phi}})$ are both non-zero elements of $\hf(-\ms,\so)$. On the other hand, both invariants are primitive elements of $\hf(-\ms,\so)$ since each contact structure is Stein-fillable. Since $\ms$ is an $L$-space, this implies that these two invariants are equal (up to sign), a contradiction. Therefore, $J=\bar{J}$, which implies that $c_1(W,J)=0$. Moreover, Ozsv{\'a}th and Szab{\'o} prove in \cite{osz2} that $b_2^+(W)=0$ for any symplectic-filling of a contact structure on an $L$-space. 

We can assume that $W$ is obtained from a 3-ball by attaching $2g$ 1-handles and $k$ 2-handles. This handle decomposition of $W$ gives rise to a surgery diagram for $\partial W = \ms$ by replacing an $n$-framed 2-handle by $n$-surgery around the corresponding attaching curve and by replacing each of the the $2g$ 1-handles with a $0$-surgery. It is an easy exercise to check (using the linking matrix for this surgery diagram) that $b_1(\ms)=0$ implies that $k\geq 2g$ and that the intersection form of $W$ is negative-definite of rank $k-2g$.

Thus, $\chi(W) = k-2g+1$, $\sigma(K) = -(k-2g)$, and Equation \ref{eqn:grading2} becomes $$d(\ms,\so) = \frac{k-2g}{4} = \frac{\chi(W)-1}{4}\geq 0.$$ \end{proof}

Theorem \ref{thm:stein} combined with Theorem \ref{thm:Grading} and Proposition \ref{prop:TE} therefore enables us to compute the Euler characteristic of any Stein-filling of $\xi_{T,\phi}$ (or, if the correction term is negative, to show that $\xi_{T,\phi}$ is not Stein-fillable) whenever $\ms$ is an $L$-space.

\newpage

\section{Applications to closed 3-braids}
\label{section:tb.3braid}
In this section, we apply our results on the $\mathbb{Q}$-graded Heegaard Floer homology of $\ms$ to knots and links with braid index at most 3. To begin with, we discuss applications of Theorem \ref{thm:Grading} to the smooth concordance order of 3-braid knots. Following this, we give a complete classification of quasi-alternating 3-braid links.

\subsection{Concordance and 3-braid knots}

\begin{definition} A knot $K \subset S^3 = \partial B^4$ is called \emph{smoothly slice} if $K$ bounds a smoothly embedded disk $D^2 \hookrightarrow B^4$.
\end{definition}

\begin{definition} 
The \emph{smooth knot concordance group}, $\mathcal{C},$ is the abelian group generated by equivalence classes of knots in $S^3$ in which $[K_1] = [K_2]$ if $K_1 \# -K_2$ is smoothly slice. The group operation is connected sum. \end{definition}

A knot $K$ is said to have concordance order $n$ if $[K]$ has order $n$ in $\mathcal{C}$. Applications of Heegaard Floer homology to the study of knot concordance began with discovery of the $\tau$ invariant by Ozsv{\'a}th and Szab{\'o}, and, independently, by Rasmussen \cite{osz10,ras}. There have since been several generalizations of their work \cite{jabnaik,grs,mo, greenejab}. 

Manolescu and Owens prove in \cite{mo} that the assignment $\delta(K) = 2d(\Sigma(K), \so)$ is a surjective group homomorphism $\delta: \mathcal{C} \rightarrow \mathbb{Z}$, where $\so$ is the unique self-conjugate $\Sc$ structure on $\Sigma(K)$. In particular, if $K$ has finite concordance order then $\delta(K) = 0$. With Theorem \ref{thm:Grading}, we have, in effect, computed $\delta(K)$ for all 3-braid knots $K\subset S^3$.

\begin{corollary}
\label{cor:conc}
Recall that in the context of 3-braids, $h = (\sigma_1\sigma_2)^3.$ Below, we list $\delta(K)$ for 3-braid knots $K$.
\begin{enumerate}
\item If $K$ is the closure of the 3-braid $h^d \cdot \x\y^{-a_1} \cdots \x\y^{-a_n}$, where the $a_i \geq 0$ and some $a_j \neq 0$, then $$\delta(K) = 
\left\{\begin{array}{ll} 
(n+4-\sum_{i=1}^n{a_i})/2 &{\text{if $d>0$ is odd}}\\
(n-4-\sum_{i=1}^n{a_i})/2 &{\text{if $d<0$ is odd}}\\
(n-\sum_{i=1}^n{a_i})/2 &{\text{if $d$ is even.}}
\end{array}\right.$$ 

\item If $K$ is the closure of the 3-braid $h^d \cdot \x^{m} \y^{-1},$ where $m = -1$ or $-3$, then $$\delta(K) = 
\left\{\begin{array}{ll} 
(m+3)/2 &{\text{if $d>0$ is odd}}\\
(m-5)/2 &{\text{if $d<0$ is odd}}\\
(m+9)/2 &{\text{if $d>0$ is even}}\\
(m+1)/2 &{\text{if $d\leq 0$ is even}}
\end{array}\right.$$ 
\end{enumerate}
\end{corollary}

\begin{remark}
The closed 3-braids specified by $h^d\cdot \x^{-1}\y^{-1}$ and $h^d\cdot \x^{-3}\y^{-1}$ are the torus knots $T(3,3d-1)$ and $T(3,3d-2)$, respectively. According to Murasugi, $\sigma(T(3,3d-1)) = 4d$ and $\sigma(T(3,3d-2)) = 4d-2$ \cite{mur}. Therefore, $T(3,3d-2)$ has infinite concordance order for all $d$, while $T(3,3d-1)$ has finite concordance order only if $d=0$. In the latter case, $T(3,3d-1)$ is the unknot, which has concordance order 1. Hence, the only 3-braid knots with unknown concordance orders are those in Corollary \ref{cor:conc}.1.

\end{remark}

Manolescu and Owens show that if $K$ is an alternating knot then $\delta(K) = -\sigma(K)/2$. Below, we see that Corollary \ref{cor:conc} provides an infinite family of knots for which $\delta(K) \neq -\sigma(K)/2$. The following computation of signatures is due to Erle \cite{erle}.

\begin{proposition}
\label{prop:sig} If $K=h^d \cdot \x\y^{-a_1} \cdots \x\y^{-a_n}$, where the $a_i \geq 0$ and some $a_j \neq 0$, then $$\sigma(K) = (-n-4d+\sum_{i=1}^n{a_i}).$$
\end{proposition}

Note that $\sigma(K)$ depends on the exponent $d$ whereas $\delta(K)$ does not. Since both are concordance invariants, a comparison of Corollary \ref{cor:conc} with Proposition \ref{prop:sig} immediately yields the following.

\begin{proposition}
\label{prop:conc2}
 If $K$ is a knot with braid index at most 3 and $K$ has finite concordance order, then $K$ is the closure of a braid $h^d \cdot \x\y^{-a_1} \cdots \x\y^{-a_n}$ where the $a_i \geq 0$, some $a_j \neq 0$, $d\in \{-1,0,1\}$, and $n+4d=\sum_{i=1}^n{a_i}.$
\end{proposition}

As we shall see in the next subsection, Proposition \ref{prop:conc2} implies that any 3-braid knot with finite concordance order is quasi-alternating. As such, it is probably the case that the $\tau$ invariant from knot Floer homology and Rasmussen's $s$ invariant from Khovanov homology (see \cite{ras3}) are both equal to $-\sigma/2$ and, hence, do not provide additional concordance information. On the other hand, when $K$ is quasi-alternating there is a fairly straightforward algorithm to compute the correction terms $d(\Sigma(K),\spc)$ for all $\Sc$ structures $\spc$ on $\Sigma(K)$ \cite{osz12}. One might hope to use these correction terms as in \cite{jabnaik,greenejab} to obtain finer concordance information about $K$ than what is provided here by $\delta$. 

\subsection{Quasi-alternating 3-braid links}
\label{ssection:tb.3braid.qa}
In this subsection, we identify all quasi-alternating 3-braid links. In order to prove that a particular 3-braid link $K$ is quasi-alternating we either note that it is alternating, or we show explicitly that it satisfies the conditions of Definition \ref{def:Q}. To prove that a link $K$ is \emph{not} quasi-alternating we either observe that $\Sigma(K)$ is not among the $L$-spaces in our classification, or we show that $\text{rk}(\kh(K)) >\text{det}(K)$. Below is the main theorem of this subsection.

\begin{theorem} 
\label{thm:QA}
The following is a complete classification of quasi-alternating links with braid index at most 3.
\begin{enumerate}
\item If $K$ is the closure of the braid $h^d \cdot \x\y^{-a_1} \cdots \x\y^{-a_n}$, where the $a_i \geq 0$ and some $a_j \neq 0$, then $K$ is quasi-alternating if and only if $d\in \{-1,0,1\}$. \\

\item If $K$ is the closure of the braid $h^d \cdot y^m,$ then $K$ is quasi-alternating if and only if either $d = 1$ and $m\in \{-1,-2,-3\}$ or $d=-1$ and $m \in \{1,2,3\}$.\\

\item If $K$ is the closure of the braid $h^d \cdot \x^m \y^{-1},$ where $m \in \{-1,-2,-3\},$ then $K$ is quasi-alternating if and only if $d\in \{0,1\}$. 

\end{enumerate}
\end{theorem}

\begin{proof}[Proof of Theorem \ref{thm:QA}] This proof is broken into three parts according to the three families of braids above.\\ 

\noindent \textbf{I.} Let $K$ be the closed 3-braid specified by $h^d \cdot \x\y^{-a_1} \cdots \x\y^{-a_n}$, where the $a_i \geq 0$ and some $a_j \neq 0$. By Theorem \ref{thm:GOFClass}, $\Sigma(K)$ is an $L$-space only if $d\in\{-1,0,1\}$. It follows that $K$ is quasi-alternating only if $d\in\{-1,0,1\}$. To see that $K$ is quasi-alternating for $d\in \{-1,0,1\}$, first observe that $K$ is alternating when $d=0$. Now, suppose $d=1$ and let $K_0$ and $K_1$ be the oriented and unoriented resolutions of $K$, respectively, at the crossing $c$ indicated in Figure \ref{fig:QARes}. Then $K_0$ is the closure of  $h \cdot \x\y^{-a_1} \cdots \x\y^{-a_n+1}$ and it is easy to see that $K_1$ is an alternating link. Moreover, $\Sigma(K_0) \cong M_{T,h \cdot xy^{-a_1} \cdots xy^{-a_n+1}}$ is obtained via $-1$-surgery around a knot in $\Sigma(K)$ and $\Sigma(K_1) \cong -Q(a_1+1,a_2,\dots,a_{n-1},a_{n-1}+1)$ is obtained via $0$-surgery around the same knot.\footnote{See \cite{osz12} for a more detailed surgery description of the branched double covers of the two resolutions of a link at a particular crossing.} Since $\text{det}(L) = |H_1(\Sigma(L);\mathbb{Z})|$, Lemma \ref{lem:H1} implies that $$\text{det}(K) = \text{det}(K_0) + \text{det}(K_1).$$ 

\begin{figure}[!htbp]
\begin{center}
\includegraphics[width=11cm]{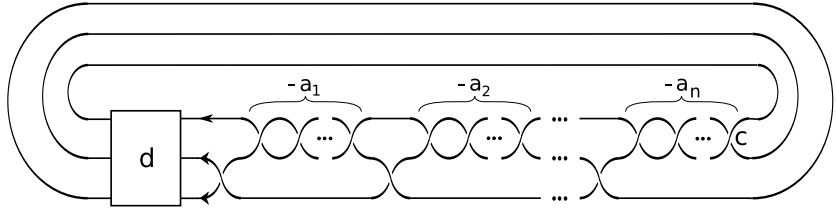}
\caption{A planar diagram for the closure of the braid $h^d \cdot \x\y^{-a_1} \cdots \x\y^{-a_n}$ with a marked crossing $c$.}
\label{fig:QARes}
\end{center}
\end{figure}

Therefore, the same induction used in the proof of Theorem \ref{thm:GOFClass} can be applied to show that $K$ is quasi-alternating. We need only to establish this fact for the closure $K'$ of the braid $h\cdot \x^n\y^{-1}$, the base case in the induction. But we have already shown that $K'$ is the torus link $T(2, 4+n)$, which is certainly alternating for $n\geq 0$. Precisely the same analysis can be carried out when $d=-1$.\\

\noindent \textbf{II.} Let $K$ be the closure of the 3-braid specified by $h^d \cdot y^m$. According to Theorem \ref{thm:GOFClass}, $\Sigma(K)$ is an $L$-space only if $d=\pm 1$. Therefore, $K$ is quasi-alternating only if $d=\pm 1$. Let $d=1$ and note that the closure of the braid $$h \cdot y^{m} = xyxyxy\cdot y^{m} = xyyxyy\cdot y^{m} = xy^2xy^{m+2}$$ has at most $1+2+1+|m+2|$ crossings. Consequently, for $m \in \{-1,-2,-3\}$, $K$ has at most 5 crossings, and is therefore alternating. The same is true when $d=-1$ and $m \in\{1,2,3\}$. 

All that remains is to show that $K$ is not quasi-alternating when $d=1$ and either $m<-3$ or $m\geq 0$ (the other cases are mirror to these). Below, we prove that $K$ is not quasi-alternating for $d=1$ and $m \geq 0$. The case $m<-3$ is virtually identical. 

Suppose that $c$ is a positive crossing in a fixed planar diagram $\Kp$ for the link $K$. We let $\Kv$ and $\Kh$ be diagrams for the oriented and unoriented resolutions of $K$ at $c$, respectively. The next proposition follows immediately from the definition of Khovanov homology, and can be found in the form below in \cite{manozs,ras2}.

\begin{proposition}
\label{prop:Exact}
If we collapse the bigrading $(i,j)$ in reduced Khovanov homology to a single grading $\delta = j-i$ then there is a long exact sequence with respect to this $\delta$ grading $$\dots \rightarrow \kh^{*-\frac{e}{2}}(\Kh) \rightarrow  \kh^{*}(\Kp) \rightarrow \kh^{*-\frac{1}{2}}(\Kv) \rightarrow \kh^{*-\frac{e}{2}-1}(\Kh) \rightarrow \dots$$ where $e$ is the difference between the number of negative crossings in the diagram $\Kh$ and those in the diagram $\Kp$.
\end{proposition}

Let $B(m)$ denote the link $h \cdot y^{m}$. Figure \ref{fig:B(m)} shows a planar diagram for $B(m)$. If we resolve this diagram at the crossing $c$, then $e=3$, and $\Kv$ and $\Kh$ are diagrams for the torus links $T(2,m+4)$ and $T(2,m)$, respectively. When $m>0$, $\kh(T(2,m)) \cong \zt^m$ (where $\zt =\zzt$ here), and is supported in $\delta$ grading $-\sigma/2 = (m-1)/2$. Hence, the long exact sequence associated to $\Kp$, $\Kv$, and $\Kh$ splits into the short exact sequences 
\begin{equation}\label{eqn:SES1}0\rightarrow \kh^{\frac{m+4}{2}}(B(m)) \rightarrow \zt^{m+4} \rightarrow \zt^m \rightarrow \kh^{\frac{m+2}{2}}(B(m))\rightarrow 0,\end{equation} and
\begin{equation}\label{eqn:SES2}0\rightarrow \kh^{\delta}(B(m)) \rightarrow 0\end{equation} 
for $\delta \neq (m+4)/2$ or $(m+2)/2$. Since the Euler characteristic of the chain complex corresponding to the short exact sequence in \ref{eqn:SES1} is zero, it follows that \begin{equation}\label{eqn:rk}\text{rk}(\kh^{\frac{m+4}{2}}(B(m))) - \text{rk}(\kh^{\frac{m+2}{2}}(B(m))) = 4.\end{equation}

\begin{figure}[!htbp]
\begin{center}
\includegraphics[width=9cm]{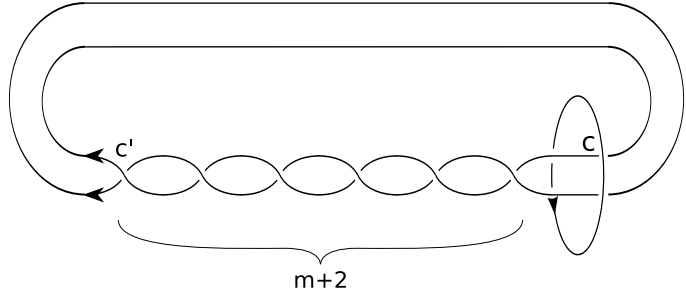}
\caption{A planar diagram for the link $B(m)$. We have indicated two positive crossings, $c$ and $c'$.}
\label{fig:B(m)}
\end{center}
\end{figure}

Let us now resolve the diagram in Figure \ref{fig:B(m)} at the crossing $c'$. In this case, $e=m+3$, and $\Kv$ and $\Kh$ are diagrams for $B(m-1)$ and the 2-component unlink $\uknot \amalg \uknot$, respectively. Recall that $\kh^{\delta}(\uknot \amalg \uknot) \cong \zt$ for $\delta = \pm 1/2$, and is zero elsewhere. From the short exact sequences above, it follows that $\kh(B(m-1))$ is supported in $\delta$-gradings $(m+3)/2$ and $(m+1)/2$ when $m>1$. Therefore, the long exact sequence associated to $\Kp$, $\Kh$, and $\Kv$ becomes \begin{eqnarray*}&0\rightarrow  &\zt \rightarrow  \kh^{\frac{m+4}{2}}(B(m))  \rightarrow  \kh^{\frac{m+3}{2}}(B(m-1))\rightarrow \\ 
&&\zt \rightarrow \kh^{\frac{m+2}{2}}(B(m)) \rightarrow  \kh^{\frac{m+1}{2}}(B(m-1)) \rightarrow 0.
\end{eqnarray*}

From this short exact sequence along with Equation \ref{eqn:rk} it follows that $$\text{rk}(\kh(B(m))) \geq \text{rk}(\kh(B(m-1))$$ when $m>1$. A computer calculation shows that $\text{rk}(\kh(B(2))) =8,$ which implies that $\text{rk}(\kh(B(m))) \geq 8$ for all $m>1$.  An addition calculation verifies that $\text{rk}(\kh(B(0))) = \text{rk}(\kh(B(1))) = 6$. On the other hand,  $\text{det}(B(m)) = 4$ for all $m$. Hence, $B(m)$ is not quasi-alternating for $m\geq 0$.\\

\noindent \textbf{III.} Let $K$ be the closed 3-braid specified by $h^d \cdot \x^m\y^{-1}$, where $m\in \{-1,-2,-3\}$. According to Theorem \ref{thm:GOFClass}, $\Sigma(K)$ is an $L$-space only if $d \in \{-1,0,1,2\}$. Thus, $K$ is quasi-alternating only if $d \in \{-1,0,1,2\}$. If $d=0$ or $1$ then $K = T(2,m)$ or $T(2,m+4)$, respectively, both of which are alternating for $m \in \{-1,-2,-3\}$. All that remains is to show that the link $K$ is not quasi-alternating when $d=-1$ or $2$. 

When $d=-1$ or $2$ and $m=-1$ or $-3$, the knots given by $h^d\cdot x^{m}y^{-1}$ are the torus knots $\pm T(3,5)$ and $\pm T(3,4)$. These torus knots are not quasi-alternating. Indeed, $\text{rk}(\kh(T(3,5)))  = 7$ while $\text{det}(T(3,5)) = 1$. Likewise, $\text{rk}(\kh(T(3,4)))  = 5$ while $\text{det}(T(3,4)) = 3$. Meanwhile, the links specified by $h^2\cdot x^{-2}y^{-1}$ and $h^{-1}\cdot x^{-2}y^{-1}$ are known as $\pm L_{9n15}$ in the link tables. These are not quasi-alternating either, as $\text{rk}(\kh(L_{9n15})) = 6$ while $\text{det}(L_{9n15}) = 2$.\\

\end{proof}

Recall that there is a spectral sequence whose $E^2$ term is $\kh(K)$ and which converges to $\hf(\Sigma(K);\zzt)$ \cite{osz12}.  This implies that $$\text{rk}(\kh(K)) \geq \text{rk}(\hf(\Sigma(K);\zzt)) \geq \text{det}(K).$$ If $\kh(K)$ is $\sigma$-thin then $\text{rk}(\kh(K)) = \text{det}(K)$, from which it follows that $\text{rk}(\hf(\Sigma(K));\zzt) = \text{det}(K)$, and, hence, that $\Sigma(K)$ is a $L$-space when working with $\zzt$ coefficients. Together with the proof of Theorem \ref{thm:QA}, this shows that the reduced Khovanov homology of a 3-braid link is $\sigma$-thin if and only if the link is quasi-alternating, as claimed in Proposition \ref{prop:thin}.

Theorem \ref{thm:QA} shows that quasi-alternating links are perhaps not as prevalent among non-alternating links as one might have thought. Indeed, for the links in Theorem \ref{thm:QA}.(2) and \ref{thm:QA}.(3), quasi-alternating is the same as alternating. Moreover, the quasi-alternating links in Theorem \ref{thm:QA}.(1) differ from closures of alternating braids by no more than a single twist in either direction. 

The following example illustrates an application of Theorem \ref{thm:QA}. In \cite{man1}, Manolescu identifies all quasi-alternating knots with 9 or fewer crossings with the exceptions of $8_{20}$ and $9_{46}$, for which $\hfk$ and $\kh$ are both $\sigma$-thin. Unfortunately, $9_{46}$ has braid index 4 and is therefore beyond the reach of Theorem \ref{thm:QA} (though, as mentioned in the introduction, Shumakovitch has shown that $9_{46}$ is not quasi-alternating). $8_{20}$, on the other hand, is a 3-braid knot. According to Bar-Natan's knot atlas, its mirror $-8_{20}$ is the closed 3-braid specified by $\y^{-3}\x\y^{3}\x$.\footnote{The knot atlas page for $8_{20}$ can be found at {\tt http://katlas.math.toronto.edu/wiki/Image:8\_20.gif}.} Up to conjugation, this braid is equivalent to 
\begin{eqnarray*}
\x\y^{3}\x\y^{-3}  &=& \x\y^{2}\cdot \y\x\y\cdot \y^{-4} = \x\y^{2}\cdot \x\y\x \cdot \y^{-4} = \x\y^{2}\x\y^{2}\cdot\y^{-1} \x \y^{-4} \\
&=& \x\y^{2}\x\y^{2}\cdot \x \y^{-5} = \x\y\x\y\x\y \cdot\x\y^{-5} = h\cdot  \x\y^{-5}. 
\end{eqnarray*}
Therefore, $-8_{20}$ is quasi-alternating according to Theorem \ref{thm:QA}.(1).

In the same way, one might hope to identify all quasi-alternating knots with 10 crossings. Of the 42 non-alternating 10-crossing knots, 32 have $\sigma$-thin knot Floer homologies, and so are potentially quasi-alternating \cite{baldgill}. Of those 32, 10 have braid index 3 and are, in fact, quasi-alternating according to our classification. See the table below for a synopsis. 

Since this manuscript first appeared, Champanerkar and Kofman have shown in \cite{ck} (incorporating work of Josh Greene) that the remaining knots in this table are quasi-alternating, with the exception of $10_{140}$, which is not quasi-alternating by forthcoming work of Shumakovitch.  \vspace{2mm}

\begin{table}[ht]
 \label{10qa}
 \begin{center}
\begin{tabular}{|c|c|c|c|} \hline
$\underline{10_{125}}$ & $10_{135}$& $10_{146}$&$\underline{10_{157}}$ \cr 
$\underline{10_{126}}$ & $10_{137}$& $10_{147}$&$10_{158}$ \cr 
$\underline{10_{127}}$ & $10_{138}$& $\underline{10_{148}}$&$\underline{10_{159}}$ \cr 
$10_{129}$ & $10_{140}$& $\underline{10_{149}}$&$10_{160}$ \cr 
$10_{130}$ & $\underline{10_{141}}$& $10_{150}$&$10_{162}$ \cr 
$10_{131}$ & $10_{142}$& $10_{151}$&$10_{163}$ \cr 
$10_{133}$ & $\underline{10_{143}}$& $\underline{10_{155}}$&$10_{164}$ \cr 
$10_{134}$ & $10_{144}$& $10_{156}$&$10_{165}$ \cr \hline

\end{tabular}
\end{center}
\vspace{3mm}
\caption{Above, we have listed the 32 $\sigma$-thin non-alternating knots with 10 crossings. Those with braid index 3 are underlined. One can show via Theorem \ref{thm:QA} that the latter are quasi-alternating.}
\end{table}

\bibliographystyle{hplain.bst}
\bibliography{References}

\end{document}